\documentclass[11pt]{amsart}
\usepackage[T1]{fontenc}
\usepackage[utf8]{inputenc}
\usepackage{lmodern,comment}
\usepackage{amsmath,amssymb,amsthm,mathtools,mathrsfs}
\usepackage[a4paper,margin=27mm]{geometry}
\usepackage{microtype,booktabs,longtable,fancyhdr,xcolor}
\definecolor{dgreen}{RGB}{0,100,50}
\usepackage[hidelinks]{hyperref}
\hypersetup{pdftitle={A C2-cofinite non-rational principal W-algebra of type D4},
 pdfauthor={Drazen Adamovic and Ivana Vukorepa},
 pdfsubject={A C2-cofinite non-rational principal W-algebra of type D4}}
\numberwithin{equation}{section}
\newtheorem{theorem}{Theorem}[section]
\newtheorem{proposition}[theorem]{Proposition}
\newtheorem{lemma}[theorem]{Lemma}
\newtheorem{corollary}[theorem]{Corollary}
\theoremstyle{definition}

\theoremstyle{remark}
\newtheorem{remark}[theorem]{Remark}
\newcommand{\CC}{\mathbb C}
\newcommand{\ZZ}{\mathbb Z}
\newcommand{\QQ}{\mathbb Q}
\newcommand{\g}{\mathfrak g}
\newcommand{\n}{\mathfrak n}
\newcommand{\h}{\mathfrak h}
\newcommand{\bfr}{\mathfrak b}
\newcommand{\nplus}{\mathfrak n_+}
\newcommand{\nminus}{\mathfrak n_-}
\newcommand{\one}{\mathbf1}
\newcommand{\W}{\mathcal W}

\newcommand{\D}{\mathbb D}
\newcommand{\HC}{\operatorname{HC}}
\newcommand{\Com}{\operatorname{Com}}
\newcommand{\Specm}{\operatorname{Specm}}
\newcommand{\End}{\operatorname{End}}

\newcommand{\wt}{\operatorname{wt}}
\newcommand{\DS}{\mathrm{DS}}
\newcommand{\reg}{\mathrm{reg}}

\newcommand{\id}{\operatorname{id}}
\allowdisplaybreaks[2]
\fancypagestyle{plain}{\fancyhf{}
 \fancyhead[R]{}
 \fancyfoot[L]{}
 \fancyfoot[R]{\thepage}}
\date{19 November 2026}
\title{A $C_2$-cofinite non-rational principal\\
 $W$-algebra of type $D_4$}
 
 \author[]{Dra\v zen  Adamovi\' c}
\address{D.A. Department of Mathematics, Faculty of Science \\
	University of Zagreb \\
	Bijeni\v cka 30 }
\email{adamovic@math.hr}

\author[]{Ivana Vukorepa}
\address{I.V. Faculty of Science \\
	University of Split \\
	Ru\dj era Bo\v{s}kovi\'ca 33}
\email{ivukorepa@pmfst.hr}

\begin{document}

\begin{abstract}
We show that a simple principal $W$-algebra can be $C_2$-cofinite
and non-rational. For $\mathcal W=\mathcal W_{-38/7}(D_4,f_{\mathrm{reg}})$
of central charge $-50$, D.~Xie predicted seventeen irreducible modules
and a logarithmic module of lowest conformal weight $-2$.
We determine the Zhu algebra $A(\mathcal W)$ and the Poisson algebra
$R_{\mathcal W}=\mathcal W/C_2(\mathcal W)$. Both have dimension $18$,
and $A(\mathcal W)\cong\CC^{16}\times\CC[\varepsilon]/(\varepsilon^2)$.
We construct all seventeen irreducible $\mathcal W$-modules
using the realization of $\mathcal W$ as a diagonal coset in
$L_{-14/3}(D_4)\otimes L_1(D_4)$ and our previous classification of
irreducible $L_{-14/3}(D_4)$-modules in the category $\mathcal O$.
We also construct an indecomposable $L_{-14/3}(D_4)$-module and use
the coset realization to obtain an indecomposable logarithmic $\mathcal W$-module.
Its lowest generalized $L(0)$-eigenspace is two-dimensional, and $L(0)$ acts
on this space by a nontrivial Jordan block with eigenvalue $-2$.
\end{abstract}

\maketitle

\section{Introduction}

The study of $C_2$-cofinite vertex operator algebras which are
not rational is motivated by logarithmic conformal field
theory. Generalized modules with a nonsemisimple action of
$L(0)$ account for the appearance of logarithmic terms in
correlation functions \cite{Gur}. At the same time,
$C_2$-cofiniteness imposes strong finiteness conditions.
For a $C_2$-cofinite vertex operator algebra $U$ of CFT type, the category
of grading-restricted generalized $U$-modules is a finite abelian
category with a braided tensor category structure \cite{Hu}.
Modular invariance also extends beyond the rational case
through the use of pseudotrace functions \cite{Mi}.

Basic examples include the triplet vertex operator algebras
$\mathcal W(p)$ and the even parts $SF(d)^+$ of symplectic
fermion vertex operator superalgebras \cite{Abe, AM08}.
A general deformation method for constructing logarithmic
modules for the triplet algebras and other subalgebras of
lattice vertex algebras was developed in \cite{AM09}.
The structure of the associative algebras $A(\mathcal W(p))$ and $A(SF(d)^+)$
was determined in \cite{AM11} and \cite{Abe,AC20}, respectively.
These works also establish the equalities
\[
 \dim A(\mathcal W(p))
 =\dim\bigl(\mathcal W(p)/C_2(\mathcal W(p))\bigr),   \quad  \dim A(SF(d)^+)
 =\dim\bigl(SF(d)^+/C_2( SF(d)^+)\bigr).
\]
These examples are closely related to the
representation theory of quantum groups at roots of unity.
 The logarithmic Kazhdan--Lusztig correspondence for the triplet
algebras $\mathcal W(p)$ was proposed in \cite{FGST1,FGST2}
and developed further in \cite{NT,CGR}.
Gannon and Negron \cite{GN} proved the corresponding braided tensor
equivalence for all $p$. Creutzig, Lentner and Rupert gave another
proof in \cite{CLR}.

In this paper we study the simple principal $W$-algebra
$\mathcal W_{-38/7}(D_4,f_{\mathrm{reg}})$.
We prove $C_2$-cofiniteness by deriving polynomial relations among the
classes of the generators in the $C_2$-algebra. To prove non-rationality,
we construct a grading-restricted generalized module for this $W$-algebra
on which $L(0)$ has a nontrivial Jordan block.

This algebra belongs to a family predicted by Xie and Yan
\cite[Table~13]{XY} to be $C_2$-cofinite. It occurs in the last
$D_N$ row of their table with $(N,n,u)=(4,4,7)$.
For this particular algebra,
Xie \cite[Section~6.1.1 and Table~7]{Xie} predicted seventeen
irreducible modules and a logarithmic module of lowest conformal weight
$-2$, using generalized affine Springer fibres.

Let $\g=\mathfrak{so}_8(\CC)$ and let $f_{\reg}$ be a principal
nilpotent element of $\g$. We write
\begin{equation}\label{eq:levels}
 \mathcal W=\W_{-38/7}(D_4,f_{\reg}).
\end{equation}
Here $\W^k(\g,f_{\reg})$ denotes the universal principal
$W$-algebra and $\W_k(\g,f_{\reg})$ denotes its simple quotient.

The coset theorems of \cite{ACK,ACL} give
\begin{equation}\label{eq:coset-intro}
 \mathcal W\cong
 \bigl(L_{-14/3}(D_4)\otimes L_1(D_4)\bigr)^{\g[t]}.
\end{equation}
The action of $\g[t]$ in this formula is diagonal.

We use the description of the maximal proper ideal of
$V^{-14/3}(D_4)$ as the ideal generated by three singular vectors
(Theorem~\ref{thm:AV-presentation}), together with the classification
of irreducible $L_{-14/3}(D_4)$-modules in category $\mathcal O$
from \cite[Proposition~3.4 and Theorem~3.6]{AV}.
Using the coset realization \eqref{eq:coset-intro}, we construct
seventeen pairwise nonisomorphic irreducible $\mathcal W$-modules
and a logarithmic $\mathcal W$-module. On its two-dimensional lowest
generalized $L(0)$-eigenspace, $L(0)$ has a nontrivial Jordan block
with eigenvalue $-2$.

Since $f_{\mathrm{reg}}$ is principal, the Zhu algebra of the universal
$W$-algebra $\mathcal W^{-38/7}(D_4,f_{\mathrm{reg}})$ is isomorphic
to the center $Z(U(\mathfrak{so}_8))$
by \cite[Theorem~4.16.3(ii)]{Ara07}.
Hence its quotient $A(\mathcal W)$ is commutative.
The actions of $A(\mathcal W)$ on these lowest conformal-weight spaces
yield a surjective algebra homomorphism
$
 A(\mathcal W)\longrightarrow
 \mathbb C^{16}\times\mathbb C[\varepsilon]/(\varepsilon^2)
$
(Theorem~\ref{thm:A-W-quotient}), and hence $\dim A(\mathcal W)\geq18$.
For any vertex operator algebra $V$, the canonical surjection
$R_V\longrightarrow\operatorname{gr}A(V)$ gives
$\dim A(V)\leq\dim R_V$ \cite[Proposition~3.3]{ALY}. Thus
\[
 18\leq\dim A(\mathcal W)\leq\dim R_{\mathcal W}.
\]

The seventeen simple $\mathcal W$-modules have lowest conformal weights
\[
 0,\quad (-1)^3,\quad -\frac97,\quad
 \left(-\frac{10}7\right)^3,\quad -\frac{11}7,\quad
 \left(-\frac{12}7\right)^3,\quad
 \left(-\frac{13}7\right)^3,\quad -2,\quad -\frac{15}7.
\]
Here a superscript $3$ denotes three modules with the same lowest
conformal weight. The characters of $A(\mathcal W)$ on the seventeen
one-dimensional lowest conformal-weight spaces are distinct, so the modules
are pairwise non-isomorphic (Theorem~\ref{thm:seventeen-simple-quotients}).
\subsection*{Main results}

We obtain the following structure theorem for
$\mathcal W=\W_{-38/7}(D_4,f_{\reg})$.

\begin{theorem}\label{thm:main}
The vertex operator algebra $\mathcal W$ is $C_2$-cofinite and
non-rational, and its central charge is $-50$. Moreover,
\[
 A(\mathcal W)\cong
 \CC^{16}\times\CC[\varepsilon]/(\varepsilon^2),
 \qquad \dim R_{\mathcal W}=18.
\]
The seventeen simple $\mathcal W$-modules in
Table~\ref{tab:seventeen-simple-modules} exhaust the irreducible
$\mathcal W$-modules and are ordinary.
There exists a grading-restricted generalized
$\mathcal W$-module with a two-dimensional lowest generalized
$L(0)$-eigenspace on which $L(0)$ acts by a nontrivial Jordan block
with eigenvalue $-2$.
\end{theorem}
Arakawa proved the rationality of the principal $W$-algebras in the
minimal series \cite{AraMin}.
Theorem~\ref{thm:main} gives a $C_2$-cofinite non-rational simple
principal $W$-algebra at $k=-38/7$. This level is non-admissible
by the Kac--Wakimoto criterion \cite[Proposition~3.1]{AraAdm},
since $k+h^\vee=4/7$ has numerator $4<h^\vee=6$.

We also give an explicit polynomial presentation of
$R_{\mathcal W}$ in
Theorem~\ref{thm:Zhu-C2-iso}. It is a graded algebra
with Hilbert series
\[
 1+q^2+3q^4+4q^6+4q^8+3q^{10}+q^{12}+q^{14}.
\]

The $C_2$-cofiniteness of $\mathcal W$ also follows from
the quasi-lisse property of $L_{-38/7}(D_4)$
(cf.~\cite[Theorem~1.1 and Appendix~B, Table~12]{ADS})
by applying quantum Hamiltonian reduction.
However, we need a more explicit description of Zhu's algebra
and the $C_2$-algebra of $\mathcal W$, so we give an independent proof.

\subsection*{Outline of the proof}

We first construct the $\mathcal W$-modules independently of
these presentations. The highest-weight vectors of the
$L_{-14/3}(D_4)$-modules from \cite{AV}, tensored with the vacuum
of $L_1(\g)$, generate ordinary $\mathcal W$-modules.
Their simple $\mathcal W$-module quotients give seventeen
distinct algebra homomorphisms $A(\mathcal W)\longrightarrow\CC$;
see Section~\ref{sec:simple-modules}.

For the construction of the logarithmic $\mathcal W$-module, put
$D=\CC[\varepsilon]/(\varepsilon^2)$ and
$\lambda_D=-\frac73\omega_2-\varepsilon\alpha_2$.
Write $A(L_{-14/3}(D_4))=U(\g)/I$. The nine Harish--Chandra
polynomials of Proposition~\ref{prop:AV-polynomials} vanish at $\lambda_D$.
It follows that $IM_D\cap D=0$ in the $\g$-module $M_D$ induced
from the generalized highest-weight space $D$.
Induction gives a $\mathbb Z_{\geq 0}$-graded weak $L_{-14/3}(D_4)$-module $E$.
The coset construction then gives
$\mathcal C=\mathcal W\cdot(D\otimes\one)\subset E\otimes L_1(\g)$,
whose lowest generalized $L(0)$-eigenspace is $D\otimes\one$.
On this space, $L(0)=-2\id+\frac37\varepsilon$.
Finite-dimensionality of its generalized conformal-weight
spaces is proved directly. Thus $\mathcal C$ is a
grading-restricted generalized $\mathcal W$-module, and its
nontrivial Jordan block proves non-rationality; see
Sections~\ref{sec:generalized} and~\ref{sec:logarithmic}.

Let $U=\W^{-38/7}(D_4,f_{\reg})$ and let $J$ be its maximal ideal.
In Section~\ref{sec:singular-relations} we derive relations in
$A(\mathcal W)$ and $R_{\mathcal W}$ from singular vectors in $J$.
Applying principal Drinfeld--Sokolov reduction to the submodules of
$V^{-38/7}(D_4)$ generated by the three singular vectors of conformal
weight $14$ from \cite[Theorem~1.2(5)]{JS}, we obtain three linearly
independent singular vectors $s_1,s_3,s_4\in J$ of conformal weight eight
and determine the action of $A(U)$ on each highest-weight line $\CC s_i$.
Using the Miura homomorphism, we give explicit formulas for the $s_i$
and express their classes $[s_i]\in A(U)$ as polynomials. These
polynomials vanish in $A(\mathcal W)$. Their components of weight eight
represent the classes $\bar s_i=s_i+C_2(U)$ in $R_U$ and give three
relations in $R_{\mathcal W}$.

The Miura homomorphism embeds $R_U$ in $R_{\mathcal H}$, where
$\mathcal H$ is the rank-four Heisenberg vertex algebra of
Section~\ref{subsec:principal-D4}. Write $q_i=b_i(-1)\one+C_2(\mathcal H)$
for the classes of its four generators. The Poisson bracket on $R_U$
vanishes, and $B(\bar a,\bar b)=\overline{a_{(1)}b}$ for $a,b\in U$
defines a biderivation on $R_U$. In the coordinates $q_i$, it is given by
\[
 B(F,G)=28\sum_{i=1}^4
 \frac{\partial F}{\partial q_i}\frac{\partial G}{\partial q_i}
 \qquad(F,G\in R_U).
\]
In Section~\ref{sec:Zhu-relations} we apply this formula with $G$ the
class of a vector in $J$. The resulting relations give
$\dim R_{\mathcal W}\le18$, hence $C_2$-cofiniteness.

The action on the two-dimensional lowest
generalized $L(0)$-eigenspace of $\mathcal C$, together with the
sixteen other algebra homomorphisms
$A(\mathcal W)\longrightarrow\CC$, gives a surjection
$A(\mathcal W)\longrightarrow\CC^{16}\times D$ by Theorem~\ref{thm:A-W-quotient}.
The dimension comparison $\dim A(\mathcal W)\leq\dim R_{\mathcal W}$
and the bound just obtained give
\[
 18\le\dim A(\mathcal W)
 =\dim\operatorname{gr}A(\mathcal W)
 \le\dim R_{\mathcal W}\le18.
\]
Equality identifies $A(\mathcal W)$ with $\CC^{16}\times D$
and proves that the polynomial relations give a presentation of
$R_{\mathcal W}$. The resulting description of $A(\mathcal W)$
also proves that the seventeen constructed simple
$\mathcal W$-modules exhaust the irreducible
$\mathcal W$-modules; see Theorem~\ref{thm:Zhu-C2-iso}.

\section{Preliminaries}\label{sec:notation}

Fix a triangular decomposition $\g=\nminus\oplus\h\oplus\nplus$
with Cartan subalgebra $\h$, and put $\bfr=\h\oplus\nplus$.
We normalize the invariant form $(\cdot\mid\cdot)$ on $\g$ so that the
roots have squared length two and use it to identify $\h$ with $\h^*$.
The dual Coxeter number is $6$. Let $W$ be the Weyl group, and put
$\rho=\tfrac12\sum_{\alpha>0}\alpha$ and
$\rho^\vee=\tfrac12\sum_{\alpha>0}\alpha^\vee$.

For a vertex operator algebra $V$, we use
\[
 Y(u,z)=\sum_{n\in\ZZ}u_{(n)}z^{-n-1},\qquad
 Y(\omega,z)=\sum_{n\in\ZZ}L(n)z^{-n-2}.
\]
For homogeneous $u\in V$, put $o(u)=u_{(\wt(u)-1)}$.
Let $O(V)$ be the span of $\sum_{j\geq0}\binom{\wt(u)}j u_{(j-2)}v$
for homogeneous $u\in V$ and $v\in V$. Write $[a]=a+O(V)$ for the
class of $a\in V$ in the Zhu algebra $A(V)=V/O(V)$. Its product is
\[
 [u]*[v]=\left[\sum_{j\geq0}\binom{\wt(u)}j u_{(j-1)}v\right]
 \qquad(u\text{ homogeneous},\ v\in V).
\]
We also write
\[
 C_2(V)=\operatorname{span}_{\CC}\{u_{(-2)}v:u,v\in V\},
 \quad R_V=V/C_2(V),\quad X_V=\Specm R_V.
\]
For $a\in V$ put $\bar a=a+C_2(V)$. The product and Poisson bracket
on $R_V$ are given by
\[
 \bar u\,\bar v=\overline{u_{(-1)}v},\qquad
 \{\bar u,\bar v\}=\overline{u_{(0)}v}.
\]
We regard $X_V$ as a reduced variety.
We say that $V$ is $C_2$-cofinite if $\dim R_V<\infty$.

\begin{lemma}\label{lem:first-product}
Let $V$ be a vertex algebra such that the Poisson bracket on $R_V$
vanishes. For $a,b\in V$, the formula
\[
 B_V(\bar a,\bar b)=\overline{a_{(1)}b}
\]
defines a symmetric biderivation of $R_V$. If $K$ is a vertex algebra
ideal of $V$ and $\bar K=(K+C_2(V))/C_2(V)$, then
\[
 B_V(R_V,\bar K)\subseteq\bar K.
\]
\end{lemma}

\begin{proof}
For every $a\in V$ we have $a_{(0)}V\subseteq C_2(V)$.
By \cite[Lemma~3.2]{ALY}, $a_{(1)}$ induces a derivation of $R_V$.
Moreover,
$\overline{a_{(1)}b}=\overline{b_{(1)}a}$
by \cite[Lemma~3.1(4)]{ALY}. Hence the product depends only on
$\bar a$ and $\bar b$ and is a derivation in each argument.
The last assertion follows from $a_{(1)}K\subseteq K$.
\end{proof}

A weak $V$-module $M$ is $\mathbb Z_{\geq 0}$-graded if it has a grading
$M=\bigoplus_{n\geq0}M(n)$ such that
$u_{(r)}M(n)\subset M(n+\wt(u)-r-1)$ for homogeneous $u\in V$,
$r\in\ZZ$ and $n\geq0$, with $M(j)=0$ for $j<0$.
The grading operator is $D_M$, defined by $D_M|_{M(n)}=n\id$.
We call $V$ rational if every $\mathbb Z_{\geq 0}$-graded weak $V$-module is completely reducible.
We say that a generalized $V$-module is grading-restricted if it is a
direct sum of finite-dimensional generalized $L(0)$-eigenspaces and
the real parts of its conformal weights are bounded below.
We call such a $V$-module logarithmic when $L(0)$ is not semisimple.

The affine Lie algebra associated with $\g$ is
\[
 \widehat{\g}=\g\otimes\CC[t,t^{-1}]\oplus\CC K.
\]
For $x\in\g$ and $r\in\ZZ$, write $x(r)=x\otimes t^r$. The Lie bracket is
\[
 [x(r),y(s)]=[x,y](r+s)+r\delta_{r+s,0}(x\mid y)K,
 \qquad [K,\widehat{\g}]=0.
\]
We say that a $\widehat{\g}$-module has level $k$ if $K$ acts as
$k\id$.
The vectors $x=x(-1)\one$ and $y=y(-1)\one$ in $V^k(\g)$ satisfy
$[x_\lambda y]=[x,y]+k(x\mid y)\lambda\one$.
For $k\neq-6$, its conformal vector is
\begin{equation}\label{eq:sugawara}
 \omega_k=\frac1{2(k+6)}\sum_{j=1}^{28}u_j(-1)u^j(-1)\one,
\end{equation}
where $\{u_j\}$ and $\{u^j\}$ are bases of $\g$ dual with respect to $(\cdot\mid\cdot)$.
The quadratic Casimir is $\Omega=\sum_{j=1}^{28}u_ju^j$.
It acts on a highest-weight vector of weight $\lambda$ by
$(\lambda\mid\lambda+2\rho)$.
The standard isomorphism $A(V^k(\g))\cong U(\g)$ \cite{FZ} identifies
$[\omega_k]$ with $\Omega/(2(k+6))$.
We write $L^k(n)=(\omega_k)_{(n+1)}$. The grading of $V^k(\g)$
is the eigenspace decomposition of $L^k(0)$.
For $\lambda\in\h^*$, let $L_k(\lambda)$ denote the simple
highest-weight $\widehat{\g}$-module of level $k$ whose highest-weight
vector has $\g$-weight $\lambda$; in particular, $L_k(0)=L_k(\g)$.
For $k\neq-6$, put
\[
 \Delta_k(\lambda)=\frac{(\lambda\mid\lambda+2\rho)}{2(k+6)},
 \qquad
 D^k_\lambda=L^k(0)-\Delta_k(\lambda)\id.
\]
The operator $D^k_\lambda$ vanishes on the highest-weight vector of
$L_k(\lambda)$ and satisfies $[D^k_\lambda,x(r)]=-r x(r)$ for
$x\in\g$, $r\in\ZZ$. Its eigenspaces give the $\ZZ_{\geq0}$-grading on
$L_k(\lambda)$.

Put $U(\g)_0=\{u\in U(\g):[h,u]=0\text{ for }h\in\h\}$.
The unshifted Harish--Chandra projection is
\begin{equation}\label{eq:HC}
 \HC:U(\g)_0\longrightarrow U(\h)=S(\h),
\end{equation}
obtained from the triangular decomposition. Thus $u\in U(\g)_0$
acts on a highest-weight vector of weight $\lambda$ by
$\HC(u)(\lambda)$.

For a finite-dimensional simple Lie algebra $\mathfrak a$ and a nilpotent
element $f\in\mathfrak a$, let $H^0_{\DS,f}$ denote the zeroth cohomology
of the Drinfeld--Sokolov complex associated with $f$. The universal
affine $W$-algebra is defined by quantum Hamiltonian reduction \cite{KRW} as
\[
 \W^k(\mathfrak a,f)=H^0_{\DS,f}\bigl(V^k(\mathfrak a)\bigr).
\]
For $\mathfrak a=\mathfrak{sl}_n$ and $f$ principal, this construction
is due to Feigin and Frenkel \cite{FF}.
We denote the simple graded quotient by $\W_k(\mathfrak a,f)$.

\subsection{The principal \texorpdfstring{$W$}{W}-algebra of type
\texorpdfstring{$D_4$}{D4}}\label{subsec:principal-D4}
Let $f_{\reg}$ be the sum of the negative simple root vectors of $\g=D_4$.
By Feigin--Frenkel duality \cite{FF91}, in the form of
\cite[Theorem~5.6]{ACL}, we identify
$\W^{-38/7}(D_4,f_{\reg})$ with $\W^{-17/4}(D_4,f_{\reg})$.
The free-field construction of the generators in type $D$ is due to
Luk'yanov and Fateev \cite{LF}. At $\alpha_0=3/\sqrt{28}$, let $R_4$
denote the weight-four generator in the normalization of Ito and Zhu
\cite[(3.16)]{IZ}. Set $Q=112R_4$.
Let $\mathcal H$ be the rank-four Heisenberg vertex algebra generated
by $b_1,\ldots,b_4$. Write $b_i(z)=Y(b_i,z)=\sum_{n\in\ZZ}b_i(n)z^{-n-1}$. Then
\[
 [b_i{}_{\lambda} b_j]=28\delta_{ij}\lambda\one .
\]
We equip $\mathcal H$ with the conformal vector $L$ in
\eqref{eq:primary-generators}.
Put $\delta_z=3\partial_z$. Define the fields $V_1,V_2,V_3,V_4$ by
\[
\begin{aligned}
 \mathcal D(z)
   &=(\delta_z+b_1(z))(\delta_z+b_2(z))
     (\delta_z+b_3(z))(\delta_z+b_4(z))\\
   &=\delta_z^4+V_1(z)\delta_z^3+V_2(z)\delta_z^2
     +V_3(z)\delta_z+V_4(z),
\end{aligned}
\]
with all powers of $\delta_z$ on the right. The constant coefficient is
\[
 V_4(z)=\mathcal D(z)1
       =(\delta_z+b_1(z))(\delta_z+b_2(z))(\delta_z+b_3(z))b_4(z).
\]

The following theorem expresses the construction of Luk'yanov and
Fateev \cite{LF} in terms of the Miura homomorphism.
Its injectivity follows from \cite[Section~5]{ACL}.
\begin{theorem} [cf. \cite{LF}, \cite{ACL}] \label{thm:miura-D4}
The Miura homomorphism
\[
 \mu:\W^{-38/7}(D_4,f_{\reg})\longrightarrow\mathcal H
\]
is injective and preserves the conformal vector. It is determined by
$Y(\mu(Q),z)=\tfrac17 V_4(z)$.
The coefficients of the poles of orders $6,4,2$ in $Y(Q,z)Y(Q,w)$
are the fields corresponding to $Q_{(5)}Q$, $Q_{(3)}Q$ and
$Q_{(1)}Q$, of conformal weights $2,4,6$, respectively. 
Together with $Q$, these vectors strongly generate
$\W^{-38/7}(D_4,f_{\reg})$ \cite[(3.16)--(3.18)]{IZ}.
\end{theorem}


We henceforth identify $\W^{-38/7}(D_4,f_{\reg})$ with its image
in $\mathcal H$.
The diagram automorphism $\sigma$ acts on
$b=(b_1,b_2,b_3,b_4)^{\mathsf T}$ by $\sigma(b)=S^{\mathsf T}b$, where
\[
 S=\frac12\begin{pmatrix}
  1& 1& 1& 1\\
  1& 1&-1&-1\\
  1&-1& 1&-1\\
 -1& 1& 1&-1
 \end{pmatrix}.
\]

\begin{proposition}\label{prop:miura-generators}
The universal $W$-algebra $\W^{-38/7}(D_4,f_{\reg})$ is strongly
generated by its conformal vector $L=\omega$ of central charge $-50$
and three primary vectors $P,Q,R$ of conformal weights $4,4,6$,
respectively. These vectors are given by
\begin{equation}\label{eq:primary-generators}
\begin{aligned}
 L&=\frac1{56}\sum_{i=1}^4 b_i(-1)^2\one
       +\frac3{28}(3b_1+2b_2+b_3)(-2)\one,\\
 Q&=\frac17(3L(-1)+b_1(-1))(3L(-1)+b_2(-1))
              (3L(-1)+b_3(-1))b_4(-1)\one,\\
 P&=-Q-2\sigma(Q),\\
 R&=P_{(1)}P+224L(-2)P-124L(-1)^2P\\
  &\quad-\frac{832}{4747}\bigl(382944L(-6)-143874L(-4)L(-2)
                        -50829L(-3)^2+14348L(-2)^3\bigr)\one .
\end{aligned}
\end{equation}
\end{proposition}
\begin{proof}
These are the generators of \cite{LF} and \cite[Section~6.3.3]{BS} at
$\alpha_0=3/\sqrt{28}$, with $P,Q,R$ normalized as
$112W_4$, $112R_4$ and $75264W_6$, respectively, in the notation of
\cite[(3.16)--(3.17), (B.10)--(B.13)]{IZ}.
With these normalizations, $Q_{(5)}Q=126464L$ and
$Q_{(3)}Q\equiv-208P$ modulo Virasoro descendants of $\one$.
Similarly, $Q_{(1)}Q\equiv\tfrac13 R$ modulo Virasoro descendants
of $\one,P,Q$.
\end{proof}

Let $U=\W^{-38/7}(\g,f_{\reg})$ and
$q_i=\overline{b_i(-1)\one}\in R_{\mathcal H}$ for $1\leq i\leq4$.
The Miura homomorphism induces an injective Poisson algebra homomorphism
\[
 \bar\mu:R_U\longrightarrow R_{\mathcal H}=\CC[q_1,\ldots,q_4],
 \qquad \bar a\longmapsto\overline{\mu(a)},
\]
with image $R_{\mathcal H}^{W(D_4)}$ \cite[(25) and Lemma~6.1(1)]{ACL}.
We identify $R_U$ with the polynomial algebra $\CC[x,y,z,t]$, where
\begin{equation}\label{eq:x-y-z-t}
\begin{aligned}
 x&=\sum_{i=1}^4q_i^2, & y&=\sum_{1\leq i<j\leq4}q_i^2q_j^2,\\
 z&=q_1q_2q_3q_4, & t&=\sum_{1\leq i<j<l\leq4}q_i^2q_j^2q_l^2.
\end{aligned}
\end{equation}
The generators $x,y,z,t$ have weights $2,4,4,6$, respectively.

For $1\leq i,j\leq4$, the $\lambda$-bracket
$[b_i{}_\lambda b_j]=28\delta_{ij}\lambda\one$ gives
$(b_i)_{(0)}b_j=0$ and $(b_i)_{(1)}b_j=28\delta_{ij}\one$.
Since the $q_i$ generate $R_{\mathcal H}$, the first identity implies that
the Poisson bracket vanishes on $R_{\mathcal H}$ and on its subalgebra $R_U$.
Hence, by Lemma~\ref{lem:first-product} the first product induces a symmetric biderivation $
B_{\mathcal H}(\bar a,\bar b)=\overline{a_{(1)}b},
\ a,b\in\mathcal H.
$ 
From the second identity we have $B_{\mathcal H}(q_i,q_j)=28\delta_{ij}$.
For $F,G\in R_{\mathcal H}$, put
$B_0(F,G)=\sum_{i=1}^4(\partial F/\partial q_i)(\partial G/\partial q_i)$.
The Leibniz rule in the second and then in the first argument gives
\[
\begin{aligned}
 B_{\mathcal H}(q_i,G)
 &=\sum_{j=1}^4 B_{\mathcal H}(q_i,q_j)\frac{\partial G}{\partial q_j}
   =28\frac{\partial G}{\partial q_i},\\
 B_{\mathcal H}(F,G)
 &=\sum_{i=1}^4\frac{\partial F}{\partial q_i}B_{\mathcal H}(q_i,G)
   =28\sum_{i=1}^4\frac{\partial F}{\partial q_i}
                    \frac{\partial G}{\partial q_i}=28B_0(F,G).
\end{aligned}
\]

For $a,b\in U$, the identity $\mu(a_{(1)}b)=\mu(a)_{(1)}\mu(b)$ gives
$\bar\mu(\overline{a_{(1)}b})=
B_{\mathcal H}(\bar\mu(\bar a),\bar\mu(\bar b))$.
Identifying $R_U$ with its image, we obtain
\begin{equation}\label{eq:first-product-miura}
 \overline{a_{(1)}b}=28B_0(\bar a,\bar b)\qquad(a,b\in U).
\end{equation}
Projecting \eqref{eq:primary-generators} to $R_U$ and using
\eqref{eq:first-product-miura} for $P_{(1)}P$, we obtain
\[
 \bar\omega=\frac{x}{56},\qquad
 \bar Q=\frac z7,\qquad
 \bar P=\frac{4y-x^2}{56},\qquad
 \bar R=\frac{92986}{1628221}x^3-\frac27xy+\frac{12}{7}t.
\]

\section{The vertex algebra \texorpdfstring{$L_{-14/3}(D_4)$}{L(-14/3)(D4)}}
\label{sec:affine-AV}
In this section we recall the structure and representation theory of the vertex algebra $L_{-14/3}(D_4)$ from our recent paper \cite{AV}.
\vskip 5mm 
We use the simple roots
\[
 \alpha_1=\epsilon_1-\epsilon_2,\qquad
 \alpha_2=\epsilon_2-\epsilon_3,\qquad
 \alpha_3=\epsilon_3-\epsilon_4,\qquad
 \alpha_4=\epsilon_3+\epsilon_4,
\]
where $(\epsilon_i\mid\epsilon_j)=\delta_{ij}$. Write
\[
 Q=\bigoplus_{i=1}^4\ZZ\alpha_i,\qquad
 Q_+=\sum_{i=1}^4\ZZ_{\geq0}\alpha_i.
\]
The fundamental weights are
\[
 \omega_1=\epsilon_1,\qquad \omega_2=\epsilon_1+\epsilon_2,
 \qquad
 \omega_3=\tfrac12(\epsilon_1+\epsilon_2+\epsilon_3-\epsilon_4),
 \qquad
 \omega_4=\tfrac12(\epsilon_1+\epsilon_2+\epsilon_3+\epsilon_4).
\]
In particular, $\rho=3\epsilon_1+2\epsilon_2+\epsilon_3$,
$(\rho\mid\rho)=14$, and $h^\vee=6$.
Write $h_i=\alpha_i^\vee$ and
$\g=\mathfrak n_-\oplus\h\oplus\mathfrak n_+$.
Let $\sigma$ be the diagram automorphism determined by
\[
 \alpha_1\longmapsto\alpha_3\longmapsto\alpha_4
 \longmapsto\alpha_1,\qquad \sigma(\alpha_2)=\alpha_2.
\]
Its Chevalley lift acts on $V^{-14/3}(\g)$.

We recall the presentation of $L_{-14/3}(\g)$ from \cite{AV}.
The normalizations of root vectors used there can be specified inside the
Clifford algebra generated by $a_i,a_i^*$, $1\leq i\leq4$, with
\[
 a_i a_j+a_j a_i=a_i^*a_j^*+a_j^*a_i^*=0,
 \qquad a_i a_j^*+a_j^*a_i=\delta_{ij}.
\]
With $:uv:=\tfrac12(uv-vu)$, take, for $i<j$,
\[
 \begin{aligned}
 e_{\epsilon_i-\epsilon_j}&=:a_i a_j^*:,&
 e_{\epsilon_i+\epsilon_j}&=:a_i a_j:,\\
 f_{\epsilon_i-\epsilon_j}&=:a_j a_i^*:,&
 f_{\epsilon_i+\epsilon_j}&=:a_j^*a_i^*:.
 \end{aligned}
\]
The Cartan elements $H_i=:a_i a_i^*:$ are dual to $\epsilon_i$.

\begin{theorem}[{\cite[Theorems 3.1 and 3.6]{AV}}]
\label{thm:AV-presentation}
There is a singular vector $v_1\in V^{-14/3}(\g)$ of conformal weight
$6$ and $\g$-weight $2\omega_1$. Set
$v_3=\sigma(v_1)$ and $v_4=\sigma^2(v_1)$. The maximal proper ideal
$J\subset V^{-14/3}(\g)$ is generated by $v_1,v_3,v_4$. Consequently,
\[
 L_{-14/3}(\g)=V^{-14/3}(\g)/\langle v_1,v_3,v_4\rangle.
\]
\end{theorem}

Here $v_1$ is taken with the normalization in the supplementary file
\texttt{D4-sing-v.nb} accompanying \cite[Theorem~3.1]{AV}. The element
of $U(\g)$ corresponding to $[v_1]$ is given in \texttt{D4-category-O.nb} accompanying
\cite[Section 3]{AV}. 

Under the canonical isomorphism $A(V^{-14/3}(\g))\cong U(\g)$,
let $r_i$ be the element of $U(\g)$ corresponding to $[v_i]$ for $i=1,3,4$. Put
\[
 R^{(i)}=U(\g)_{\mathrm{ad}}r_i\quad(i=1,3,4),
 \qquad R=R^{(1)}\oplus R^{(3)}\oplus R^{(4)}.
\]
For $\mu\in\h^*$, write $R_\mu=\{r\in R:[h,r]=\mu(h)r
\text{ for all }h\in\h\}$ and $R^{(i)}_\mu=R^{(i)}\cap R_\mu$.
The three summands are isomorphic to the finite-dimensional simple
$\g$-modules of highest weights $2\omega_1,2\omega_3,2\omega_4$,
respectively. In particular, $\dim R^{(i)}_0=3$ for $i=1,3,4$.
Theorem~\ref{thm:AV-presentation} and the induced homomorphism
$A(V^{-14/3}(\g))\longrightarrow A(L_{-14/3}(\g))$ give
\begin{equation}\label{eq:A-affine-presentation}
 A(L_{-14/3}(\g))\cong U(\g)/I,
 \qquad I=U(\g)R.
\end{equation}
For a root vector $f$, write $\operatorname{ad}(f)$ for its adjoint
operator on $U(\g)$.

\begin{proposition}[{\cite[Lemma 3.3]{AV}}]\label{prop:AV-polynomials}
Let $p_1,p_2,p_3\in S(\h)$ be the polynomials specified in
\cite[Lemma~3.3]{AV}. They form a basis of $\HC(R^{(1)}_0)$.
Moreover,
\[
 \operatorname{HC}(R_0)=
 \operatorname{span}_{\mathbb C}
 \{\sigma^j(p_i):1\leq i\leq3,\ 0\leq j\leq2\}.
\]
For $1\leq i\leq3$ and $j=1,2$, set $p_{3j+i}=\sigma^j(p_i)$.
\end{proposition}

For $\lambda\in\h^*$, the simple highest weight
$\g$-module with highest weight $\lambda$ is an
$A(L_{-14/3}(\g))$-module, under \eqref{eq:A-affine-presentation}, if and only if
\[
 \sigma^j(p_i)(\lambda)=0
 \qquad(1\leq i\leq3,\ 0\leq j\leq2).
\]
The zero set consists of the $405$ weights classified in
\cite[Proposition 3.4 and Theorem 3.6]{AV}. 
The following lemma will be used in Section \ref{sec:generalized}.

\begin{lemma}\label{lem:AV-tangent}
Let $t$ be an indeterminate and put
\[
 \lambda(t)=-\tfrac73\omega_2-t\alpha_2
 =t\omega_1+(-\tfrac73-2t)\omega_2+t\omega_3+t\omega_4.
\]
Then
\begin{align*}
 p_1(\lambda(t))&=4t^2(3t+2)^2(3t+4)^2,\\
 p_2(\lambda(t))&=-4t^2(3t+2)^2(3t+4)^2,\\
 p_3(\lambda(t))&=0.
\end{align*}
Consequently, $p_i(\lambda(\varepsilon))=0$ in
$\CC[\varepsilon]/(\varepsilon^2)$ for $1\leq i\leq9$.

\end{lemma}

\begin{proof}
Substituting $h_1=h_3=h_4=t$ and $h_2=-7/3-2t$ into the polynomials
of \cite[Lemma~3.3]{AV} gives the three identities. Since
$\sigma\lambda(t)=\lambda(t)$, the same identities apply to the
triality translates. Setting $t=\varepsilon$ proves the last assertion.
\end{proof}

\section{Coset realization of \texorpdfstring{$\W_{-38/7}(D_4,f_{\reg})$}{W(-38/7)(D4, f reg)}}
\label{sec:coset}
We realize the simple principal $W$-algebra $\W_{-38/7}(D_4,f_{\reg})$ as a commutant.
\vskip 5mm

The diagonal action of $\widehat{\g}$ defines a vertex algebra
homomorphism
\[
 \varphi:V^{-11/3}(\g)\longrightarrow
 L_{-14/3}(\g)\otimes L_1(\g),\qquad
 x(-1)\one\longmapsto x(-1)\one\otimes\one+
 \one\otimes x(-1)\one.
\]
Write $x^\varphi(n)=\varphi(x(-1)\one)_{(n)}$ for $x\in\g$, $n\in\ZZ$.
The commutant is
\begin{equation}\label{eq:coset-definition}
 \begin{aligned}
 &\Com\bigl(\varphi(V^{-11/3}(\g)),
 L_{-14/3}(\g)\otimes L_1(\g)\bigr)\\
 &\qquad=\{v\in L_{-14/3}(\g)\otimes L_1(\g):
x^\varphi(n)v=0\text{ for all }x\in\g,\ n\geq0\}.
 \end{aligned}
\end{equation}
With $\omega_k$ the Sugawara vector in \eqref{eq:sugawara}, the conformal vector of the commutant is
\begin{equation}\label{eq:coset-conformal}
 \omega^{\mathcal W}=\omega_{-14/3}\otimes\one+
 \one\otimes\omega_1-\varphi(\omega_{-11/3}).
\end{equation}

By \cite[Theorem~4.3]{ACK}, the coset is isomorphic to
$\W_{-17/4}(\g,f_{\reg})$. Feigin--Frenkel duality identifies this
algebra with $\W_{-38/7}(\g,f_{\reg})$. We recall the argument
of \cite{ACK} in our notation.

\begin{proposition}\label{prop:coset}
There is an isomorphism of conformal vertex algebras
\[
 \Com\bigl(\varphi(V^{-11/3}(\g)),
 L_{-14/3}(\g)\otimes L_1(\g)\bigr)
 \cong\W_{-38/7}(\g,f_{\reg}).
\]
In particular, $\mathcal W$ has central charge $-50$.
\end{proposition}
\begin{proof}
Since $-14/3+6=4/3\notin\QQ_{\leq0}$,
\cite[Theorem 8.7(2)]{ACL} gives
\[
 \W^{-38/7}(\g,f_{\reg})\cong
 (V^{-14/3}(\g)\otimes L_1(\g))^{\g[t]}.
\]
The diagonal action of $\widehat{\g}$ has level $k_0=-11/3$.
Let $KL^{k_0}$ denote the category of smooth $\widehat{\g}$-modules of
level $k_0$ on which $\g[t]$ acts locally finitely. Both tensor products
$V^{-14/3}(\g)\otimes L_1(\g)$ and
$L_{-14/3}(\g)\otimes L_1(\g)$ belong to this category. Since $k_0+6=7/3\notin\QQ_{\leq0}$, the
$\widehat{\g}$-module $V^{k_0}(\g)$ is projective in $KL^{k_0}$
\cite[Lemma 2.1]{ACL}. Applying
$\operatorname{Hom}_{\widehat\g}(V^{k_0}(\g),-)$ to the natural
surjection onto $L_{-14/3}(\g)\otimes L_1(\g)$ therefore gives a surjection
\begin{equation}\label{eq:universal-to-coset}
 \W^{-38/7}(\g,f_{\reg})\longrightarrow
 (L_{-14/3}(\g)\otimes L_1(\g))^{\g[t]}.
\end{equation}

The vertex algebra $L_{-14/3}(\g)\otimes L_1(\g)$ is $\ZZ_{\geq0}$-graded by
$L^{-14/3}(0)\otimes\id+\id\otimes L^1(0)$, its degree-zero subspace
is $\CC\one$, and its homogeneous subspaces are finite-dimensional.
It has a nondegenerate invariant bilinear form.
The diagonal vectors $\varphi(x(-1)\one)$, for $x\in\g$,
are primary of conformal weight one.
The simplicity theorem \cite[Theorem 4.1]{ACK}, applied at level
$-11/3$, shows that
$(L_{-14/3}(\g)\otimes L_1(\g))^{\g[t]}$ is simple.
Thus \eqref{eq:universal-to-coset} identifies this commutant with
the simple quotient of $\W^{-38/7}(\g,f_{\reg})$.
The conformal vectors are identified by the coset theorem.
By \eqref{eq:coset-conformal}, the central charge is $c_{\mathcal W}=-98+4+44=-50$.
\end{proof}

\section{A generalized highest-weight module for \texorpdfstring{$L_{-14/3}(D_4)$}{L(-14/3)(D4)}}
\label{sec:generalized}
{We construct the generalized highest-weight $L_{-14/3}(D_4)$-module used
in Section~\ref{sec:logarithmic} to obtain a logarithmic $\mathcal W$-module.
Let $D=\CC[\varepsilon]/(\varepsilon^2)$ and set
$\lambda_0=-\tfrac73\omega_2$ and $\lambda_D=\lambda_0-\varepsilon\alpha_2$.
The weight $\lambda_0$ occurs in the classification of \cite{AV}.
Regard $D$ as a $\bfr$-module by setting
\[
 \n_+D=0,\qquad h\cdot d=\lambda_D(h)d
 \quad(h\in\h,\ d\in D).
\]
Thus $h\cdot\varepsilon=\lambda_0(h)\varepsilon$ and
$h\cdot 1=\lambda_0(h)-\varepsilon\alpha_2(h)$.
In particular, $\varepsilon$ is a highest-weight vector of weight
$\lambda_0$, and $1$ is a generalized highest-weight vector.

Consider the induced $\g$-module and its quotient
\[
 M_D=U(\g)\otimes_{U(\bfr)}D\cong U(\n_-)\otimes D,
 \qquad T=M_D/IM_D,
\]
where $I$ is the ideal in \eqref{eq:A-affine-presentation}.
Then $T$ is an $A(L_{-14/3}(\g))$-module. We shall show that the
canonical map $D\longrightarrow T$ is injective.
	\begin{lemma}\label{lem:R0}
We have $IM_D\cap D=R_0D$, where $R_0$ is the zero-weight space for the adjoint action of $\h$ on $R$.
	\end{lemma}
    \begin{proof}
Since $I=U(\g)R$ is a two-sided ideal and $D$ generates $M_D$,
we have $IM_D=U(\g)RD$. The space $R$ is stable under the
adjoint action of $\g$, and $D$ is a $\bfr$-module. Thus $RD$ is
stable under $\bfr$, since $x(rd)=[x,r]d+r(xd)$ for $x\in\bfr$,
$r\in R$ and $d\in D$. It follows that $IM_D=U(\nminus)RD$.

As $M_D\cong U(\nminus)\otimes D$, its generalized $\h$-weight
spaces are finite-dimensional, its weights lie in $\lambda_0-Q_+$,
and its generalized weight space of weight $\lambda_0$ is $D$.
For the adjoint weight decomposition $R=\bigoplus_{\mu\in Q}R_\mu$,
we have $R_\mu D\subseteq (M_D)_{\lambda_0+\mu}^{\mathrm{gen}}$.
Since $\nminus$ lowers weights, the generalized weight space of
weight $\lambda_0$ in $U(\nminus)RD$ is precisely $R_0D$.
Hence $IM_D\cap D=R_0D$.
    \end{proof}
    \begin{proposition}\label{prop-inj-equiv}
    	The following statements are equivalent:
    	\begin{enumerate}
    		\item The canonical map $D \longrightarrow T$ is injective.
    		\item $p_i (\lambda_D) = 0$, for all $i=1,\ldots,9$.
    	\end{enumerate}
    \end{proposition}
\begin{proof}
The canonical map $D\longrightarrow T$ has kernel $IM_D\cap D=R_0D$
by Lemma~\ref{lem:R0}. Since $\nplus D=0$, the Harish--Chandra
projection gives
\[
 r\cdot d=\HC(r)(\lambda_D)d \qquad (r\in R_0,\ d\in D).
\]
By \cite[Lemma~3.3]{AV}, the polynomials $p_1,\ldots,p_9$ span
$\HC(R_0)$. Thus $R_0D=0$ if and only if all nine polynomials
vanish at $\lambda_D$.
\end{proof}
}

\begin{proposition}\label{prop:T}
The canonical map $D\longrightarrow T$ is injective.
Thus $T$ is a nonzero $A(L_{-14/3}(\g))$-module containing the 
two-dimensional generalized highest-weight space $D$. Its generalized $\h$-weight
spaces are finite-dimensional, and multiplication by $\varepsilon$
defines a nonzero square-zero endomorphism of $T$.
Moreover,
\begin{equation}\label{eq:Casimir-D}
 \Omega|_T=-\frac{112}{9}\id+\frac83\varepsilon.
\end{equation}
\end{proposition}
\begin{proof}
By Lemma~\ref{lem:AV-tangent}, $p_i(\lambda_D)=0$ in $D$ for $1\leq i\leq9$.
Hence $D\longrightarrow T$ is injective by Proposition~\ref{prop-inj-equiv}.
Multiplication by $\varepsilon$ on $M_D$ commutes with $U(\g)$ and
preserves $IM_D$. It therefore induces an endomorphism of $T$ whose
restriction to $D$ sends $1$ to $\varepsilon\ne0$. The generalized
$\h$-weight spaces of $T$ are finite-dimensional because those of $M_D$ are.

For $\lambda(t)=\lambda_0-t\alpha_2$, expansion of the bilinear expression
using $\lambda_0=-7\omega_2/3$ and $(\alpha_2\mid\alpha_2)=2$ gives
\[
 (\lambda(t)\mid\lambda(t)+2\rho)
 =-\frac{112}{9}+\frac83t+2t^2.
\]
Thus \eqref{eq:Casimir-D} holds on the generalized highest-weight space.
Both $\Omega$ and multiplication by $\varepsilon$ commute with
$U(\g)$, and that space generates the $\g$-module $T$, so the identity holds on $T$.
\end{proof}

\begin{proposition}\label{prop:affine-module}
There exists an indecomposable $\mathbb Z_{\geq 0}$-graded weak
$L_{-14/3}(\g)$-module $E=\bigoplus_{n\geq0}E(n)$, generated by
$E(0)=T$, on which multiplication by $\varepsilon$ extends to a
square-zero endomorphism.
For every $n\geq0$, it satisfies
\begin{equation}\label{eq:E-L0}
 L^{-14/3}(0)|_{E(n)}=\left(n-\frac{14}{3}\right)\id+\varepsilon.
\end{equation}
Its grading operator is
$
 D_E=L^{-14/3}(0)+\frac{14}{3}\id-\varepsilon.
$
Each generalized $\h$-weight space of $E(n)$ is finite-dimensional.
\end{proposition}
\begin{proof}
By Zhu's theorem \cite[Theorem~2.2.1]{Zhu}, there exists a
$\mathbb Z_{\geq 0}$-graded weak $L_{-14/3}(\g)$-module whose
degree-zero subspace is $T$ as an $A(L_{-14/3}(\g))$-module.
Let $E$ be the submodule generated by $T$, and let $D_E$ act
as $n\,\id$ on $E(n)$.

For $x\in\g$ and $r\in\ZZ$, we have
$[L^{-14/3}(0),x(r)]=[D_E,x(r)]=-r x(r)$. Thus
\[
 N=L^{-14/3}(0)-D_E+\frac{14}{3}\id
\]
commutes with every $x(r)$ and hence with the $L_{-14/3}(\g)$-action.
By \eqref{eq:sugawara} and
\eqref{eq:Casimir-D}, its restriction to $T$ is multiplication
by $\varepsilon$. Thus $N^2$ vanishes on $T$, and hence on $E$,
since $T$ generates $E$. Writing $\varepsilon$ for $N$ gives
\eqref{eq:E-L0}.

Since $L_{-14/3}(\g)$ is generated by $\g$, there is a
surjection of graded vector spaces
\begin{equation}\label{eq:E-spanning}
 U(\g\otimes t^{-1}\CC[t^{-1}])\otimes T\longrightarrow E,
\end{equation}
where $\deg x(-r)=r$ for $r\geq1$.
The fixed-degree subspaces of the first tensor factor in
\eqref{eq:E-spanning} are finite-dimensional.
Each such subspace has only finitely many $\h$-weights. Since the
generalized $\h$-weight spaces of $T$ are finite-dimensional,
\eqref{eq:E-spanning} gives $\dim E(n)^{\mathrm{gen}}_\mu<\infty$
for every $n\geq0$ and $\mu\in\h^*$.

Every $L_{-14/3}(\g)$-module endomorphism of $E$ preserves
$E(0)$, the generalized $L^{-14/3}(0)$-eigenspace of eigenvalue
$-14/3$, and its generalized $\h$-weight space $D$ of weight
$\lambda_0$. Since $h_1$ acts on $D$ as multiplication by
$\varepsilon$, restriction to $D$ defines an algebra homomorphism
\[
 \operatorname{End}_{L_{-14/3}(\g)}(E)
 \longrightarrow \operatorname{End}_{\CC[\varepsilon]}(D)\cong D.
\]
This homomorphism is injective because $D$ generates $E$.
The algebra $D$ has no nontrivial idempotents, so neither does
$\operatorname{End}_{L_{-14/3}(\g)}(E)$. Thus $E$ is indecomposable.
\end{proof}

\section{A logarithmic \texorpdfstring{$\W_{-38/7}(D_4,f_{\reg})$}{W(-38/7)(D4, f reg)}-module}
\label{sec:logarithmic}

We construct a logarithmic $\W_{-38/7}(D_4,f_{\reg})$-module
of lowest conformal weight $-2$, as predicted in \cite{Xie}.
\vskip 5mm

Using Proposition~\ref{prop:coset}, regard $\mathcal W$ as the
coset subalgebra of $L_{-14/3}(\g)\otimes L_1(\g)$.
Then $E\otimes L_1(\g)$ is a weak
$L_{-14/3}(\g)\otimes L_1(\g)$-module and, by restriction,
a weak $\mathcal W$-module.
Put $v=1\otimes\one$ and $w=\varepsilon\otimes\one$ in
$D\otimes\one\subset E\otimes L_1(\g)$. Define
\begin{equation}\label{eq:C}
 \mathcal C=\mathcal W\cdot( D\otimes\one)
 \subset E\otimes L_1(\g).
\end{equation}
Multiplication by $\varepsilon$ on the first tensor factor preserves
$\mathcal C$ and commutes with the $\mathcal W$-action. Put
$
 D_{\mathrm{tot}}=D_E\otimes\id+\id\otimes L^1(0)$. 
This operator defines a $\ZZ_{\geq0}$-grading on $E\otimes L_1(\g)$.

\begin{theorem}\label{thm:log}
The space $\mathcal C$ is a grading-restricted generalized $\mathcal W$-module.
For $n\geq0$, put $\mathcal C(n)=\{c\in\mathcal C:D_{\mathrm{tot}}c=nc\}$.
Then
\[
 \mathcal C=\bigoplus_{n\geq0}\mathcal C(n),\qquad
 \mathcal C(0)=\CC v\oplus\CC w,
\]
and
\begin{equation}\label{eq:C-L0}
 L^{\mathcal W}(0)|_{\mathcal C(n)}
 =(n-2)\id+\frac37\varepsilon.
\end{equation}
In particular,
\begin{equation}\label{eq:Jordan}
 (L^{\mathcal W}(0)+2)v=\frac37w,
 \qquad (L^{\mathcal W}(0)+2)w=0,
 \qquad w\neq0.
\end{equation}
Moreover, $\End_{\mathcal W}(\mathcal C)\cong D$, so $\mathcal C$
is indecomposable.
\end{theorem}
\begin{proof}
For $x\in\g$ and $r\in\ZZ$, the diagonal operator
$x^\varphi(r)$ acts on $E\otimes L_1(\g)$ by
\[
 x^\varphi(r)(e\otimes a)
 =x(r)e\otimes a+e\otimes x(r)a.
\]
Since $\mathcal W$ is the commutant of $\varphi(V^{-11/3}(\g))$,
each $x^\varphi(r)$ commutes with all modes of $\mathcal W$.
The action on $D\otimes\one$ therefore gives, on the
$\mathcal W$-module $\mathcal C$ generated by this subspace,
\begin{equation}\label{eq:C-highest}
\begin{aligned}
 x^\varphi(r)\mathcal C&=0
    &&(x\in\g,\ r\geq1),\\
 x^\varphi(0)\mathcal C&=0
    &&(x\in\nplus),\\
 h^\varphi(0)|_{\mathcal C}
    &=\lambda_0(h)\id-\alpha_2(h)\varepsilon
    &&(h\in\h).
\end{aligned}
\end{equation}
Here $\varepsilon$ acts on the first tensor factor.

On $D\otimes\one$, the operators $x^\varphi(0)$ act as
$x(0)$ on $D$, since $x(0)\one=0$ in $L_1(\g)$.
Thus the Sugawara formula at level $-11/3$ and
\eqref{eq:Casimir-D} give
\[
 \left.\varphi(\omega_{-11/3})_{(1)}\right|_{D\otimes\one}
 =\frac3{14}\left(-\frac{112}{9}\id+\frac83\varepsilon\right)
 =-\frac83\id+\frac47\varepsilon.
\]
Both $\varphi(\omega_{-11/3})_{(1)}$ and $\varepsilon$ commute
with all modes of $\mathcal W$, so the same identity holds on
$\mathcal C$. By Proposition~\ref{prop:affine-module},
\[
 L^{-14/3}(0)\otimes\id+\id\otimes L^1(0)
 =D_{\mathrm{tot}}-\frac{14}{3}\id+\varepsilon.
\]
It follows from \eqref{eq:coset-conformal} that on $\mathcal C$,
$
 L^{\mathcal W}(0)=D_{\mathrm{tot}}-2\id+\frac37\varepsilon$. 
In particular, $D_{\mathrm{tot}}$ preserves $\mathcal C$.
Restriction to $\mathcal C(n)$ gives \eqref{eq:C-L0}.

By \eqref{eq:C-highest}, $\mathcal C$ lies in the generalized
$\h$-weight space of weight $\lambda_0$ for the operators
$h^\varphi(0)$. In total degree $n$, the tensor product is
\[
 \bigoplus_{j=0}^n E(j)\otimes L_1(\g)_{n-j},
\]
where $L_1(\g)_m$ denotes the conformal-weight $m$ subspace.
Each $L_1(\g)_m$ is finite-dimensional and has only finitely
many $\h$-weights. For each such weight $\mu$, the generalized
$\h$-weight space of weight $\lambda_0-\mu$ in $E(j)$ is
finite-dimensional by Proposition~\ref{prop:affine-module}.
Hence $\dim\mathcal C(n)<\infty$.

In total degree zero the tensor product is $T\otimes\one$,
whose generalized $\h$-weight space of weight $\lambda_0$ is
$D\otimes\one$. Since $D\otimes\one\subset\mathcal C$, we have
$\mathcal C(0)=D\otimes\one=\CC v\oplus\CC w$.
Now \eqref{eq:Jordan} follows from \eqref{eq:C-L0} and
$\varepsilon v=w\ne0$, $\varepsilon w=0$.
The spaces $\mathcal C(n)$ are the generalized conformal-weight
spaces of weights $n-2$, so $\mathcal C$ is grading-restricted.

Every $\mathcal W$-module endomorphism of $\mathcal C$ preserves
$\mathcal C(0)=D\otimes\one$ and commutes there with $\varepsilon$
by \eqref{eq:C-L0}. Restriction to $\mathcal C(0)$ gives an
injective homomorphism $\End_{\mathcal W}(\mathcal C)\to\End_D(D)=D$,
since $\mathcal C(0)$ generates $\mathcal C$. Its image contains
multiplication by $\varepsilon$, which acts on all of $\mathcal C$.
Thus $\End_{\mathcal W}(\mathcal C)\cong D$. Since $D$ has no
nontrivial idempotents, $\mathcal C$ is indecomposable.
\end{proof}

\begin{corollary}\label{cor:nonrational}
The vertex operator algebra $\mathcal W$ is not rational.
\end{corollary}
\begin{proof}
The weak $\mathcal W$-module $\mathcal C$ is $\mathbb Z_{\geq 0}$-graded and $L(0)$ is not semisimple
on $\mathcal C(0)$.
If a vertex operator algebra $U$ is rational, every irreducible $\mathbb Z_{\geq 0}$-graded weak
$U$--module is ordinary \cite[Theorem 8.1]{DLM}, and every $\mathbb Z_{\geq 0}$-graded weak
$U$--module is a direct sum of irreducible $U$--modules.
This contradicts \eqref{eq:Jordan}.
\end{proof}

\begin{remark}
As a $\widehat{\g}$--module, $E$ is generated by the generalized
highest-weight space $D$. It does not belong to the usual category
$\mathcal O$ for $\widehat{\g}$, since
$\h$ does not act semisimply on $D$.
In the basis $1,\varepsilon$ of $D$, for example,
\[
 h_i\cdot1=\varepsilon\quad(i=1,3,4),\qquad
 h_2\cdot1=-\frac73\,1-2\varepsilon.
\]
The space $E(0)=T$ is infinite-dimensional; the subspace used in
\eqref{eq:C} is its two-dimensional generalized highest-weight space.
\end{remark}

\section{Simple \texorpdfstring{$\W_{-38/7}(D_4,f_{\reg})$}{W(-38/7)(D4, f reg)}-modules}
\label{sec:simple-modules}
In this section we construct seventeen pairwise non-isomorphic
simple $\mathcal W$-modules from $L_{-14/3}(D_4)$-modules using
the coset realization. Their completeness will follow from the
determination of $A(\mathcal W)$ in
Section~\ref{sec:singular-relations}.
\vskip 5mm

Let $\Lambda\subset\h^*$ be the set of $405$ highest weights
in the classification of irreducible $L_{-14/3}(\g)$-modules in
category $\mathcal O$; see
\cite[Proposition~3.4 and Theorem~3.6]{AV}.
For $\lambda\in\Lambda$, let $v_\lambda$ be a highest-weight vector
of $L_{-14/3}(\lambda)$ and set
\[
  \mathcal C_\lambda
  =\mathcal W\cdot(v_\lambda\otimes\mathbf1)
  \subset L_{-14/3}(\lambda)\otimes L_1(\g).
\]
The action is given by the coset realization. On
$L_{-14/3}(\lambda)\otimes L_1(\g)$ we use the grading operator
\[
 D_{\lambda,\mathrm{tot}}
 =D^{-14/3}_\lambda\otimes\id+\id\otimes L^1(0).
\]

\begin{proposition}\label{prop:cyclic-coset-modules}
For $\lambda\in\Lambda$, the $\mathcal W$-module $\mathcal C_\lambda$ is ordinary
with one-dimensional lowest conformal-weight space and lowest
conformal weight
\begin{equation}\label{eq:coset-lowest-weight}
 h_\lambda=\frac9{56}(\lambda\mid\lambda+2\rho).
\end{equation}
It has a unique irreducible $\mathcal W$-module quotient, denoted by $S_\lambda$.
For $\lambda,\mu\in\Lambda$,
\begin{equation}\label{eq:simple-coset-equivalence}
 S_\lambda\cong S_\mu
 \quad\Longleftrightarrow\quad
 \lambda+\rho\in W(\mu+\rho).
\end{equation}
\end{proposition}

\begin{proof}
Put $u_\lambda=v_\lambda\otimes\one$. On
$L_{-14/3}(\lambda)\otimes L_1(\g)$, the tensor product action gives
\[
 x^\varphi(r)u_\lambda=x(r)v_\lambda\otimes\one+
 v_\lambda\otimes x(r)\one
 \qquad(x\in\g,\ r\in\ZZ).
\]
Since $\mathcal W$ is the commutant of $\varphi(V^{-11/3}(\g))$,
the operators $x^\varphi(r)$ commute with the
$\mathcal W$-action. The highest-weight conditions on $u_\lambda$
therefore give
\[
 x^\varphi(r)c=0\quad(x\in\g,\ r\geq1),\qquad
 x^\varphi(0)c=0\quad(x\in\nplus),\qquad
 h^\varphi(0)c=\lambda(h)c\quad(h\in\h)
\]
for every $c\in\mathcal C_\lambda$.
It follows from the Sugawara formula that
$\varphi(\omega_{-11/3})_{(1)}$ acts on $\mathcal C_\lambda$
by $\Delta_{-11/3}(\lambda)$. Hence \eqref{eq:coset-conformal} gives
\[
 L^{\mathcal W}(0)
 =D_{\lambda,\mathrm{tot}}
  +\bigl(\Delta_{-14/3}(\lambda)-\Delta_{-11/3}(\lambda)\bigr)\id
 =D_{\lambda,\mathrm{tot}}
  +\frac9{56}(\lambda\mid\lambda+2\rho)\id
 \qquad\text{on }\mathcal C_\lambda.
\]

In each total degree, the $\h$-weight space of weight $\lambda$
in the tensor product is finite-dimensional: the graded components
of $L_1(\g)$ are finite-dimensional, and those of
$L_{-14/3}(\lambda)$ have finite-dimensional $\h$-weight spaces.
In degree zero this weight space is precisely $\CC u_\lambda$.
Thus $\mathcal C_\lambda$ is ordinary, with lowest conformal weight
$h_\lambda$ given by \eqref{eq:coset-lowest-weight}.
Since $u_\lambda$ generates $\mathcal C_\lambda$, every proper
submodule has zero component of degree zero. Their sum is therefore
the unique maximal submodule.

It remains to determine the character of $A(\mathcal W)$ on
$\CC u_\lambda$. We first work over the coefficient ring
$\mathscr R$ of \cite[Definition~2.4]{ACL}, with level parameter $K$,
and put $\kappa=(K+6)/(K+7)$. The coset inclusion of
\cite[Theorem~8.7]{ACL} preserves conformal weights, since the
Sugawara operator $L(0)$ for the diagonal action vanishes on the
commutant. The induced homomorphism on Zhu algebras, followed by
evaluation of $A(L_1(\g))$ on its vacuum line, gives
\[
 A(\W^{\kappa-6}(\g,f_{\reg}))
 \longrightarrow U(\g)\otimes\mathscr R.
\]
Its image lies in $Z(U(\g))\otimes\mathscr R$. Indeed, the coset
commutes with the diagonal $\g$-action. After evaluation on the vacuum
line, the commutator with $x(0)$ is the commutator with $x$ in $U(\g)$
for every $x\in\g$. By the Harish--Chandra
isomorphism, each element of the image acts on a highest-weight
line of weight $\lambda$ by a polynomial in $\lambda$ with
coefficients in $\mathscr R$.

Let $\Upsilon_Z:A(\W^{\kappa-6}(\g,f_{\reg}))\to S(\h)\otimes\mathscr R$
be the homomorphism induced by the Miura homomorphism. For
$\nu\in\h^*\otimes\mathscr R$, put $\chi_\nu=\operatorname{ev}_\nu\circ\Upsilon_Z$
as in \cite[(27)]{ACL}. For irrational $K$ and dominant integral
$\lambda$, the branching formula \cite[Theorem~11.1]{ACL}, with highest
weight $\lambda$ at levels $K$ and $K+1$ and the vacuum module at level
one, identifies the resulting character on the highest-weight line with
$\chi_{(1-\kappa)\lambda}$.
For each element of $A(\W^{\kappa-6}(\g,f_{\reg}))$ over $\mathscr R$,
both character values are polynomial in $\lambda$ with coefficients
in $\mathscr R$. At each irrational $K$, their difference vanishes
on the dominant integral weights, which are Zariski dense in $\h^*$.
Its coefficients therefore vanish for every irrational $K$ and hence
vanish in $\mathscr R$. The elements of $\mathscr R$ are regular
at $K=-14/3$ since $K+6=4/3\notin\QQ_{\leq0}$ at this level.
After specialization and passage to
$L_{-14/3}(\lambda)\otimes L_1(\g)$, the character of the
algebra $A(\W^{-38/7}(\g,f_{\reg}))$ on $\CC u_\lambda$ is
\begin{equation}\label{eq:A-coset-character}
 \chi_{3\lambda/7}.
\end{equation}

By the proof of \cite[Proposition~6.2]{ACL}, this is the central
character of the $\g$-Verma module of highest weight
$\frac37\lambda-\frac47\rho$. Its shifted Harish--Chandra parameter is
\[
 \frac37\lambda+\rho-\frac47\rho=\frac37(\lambda+\rho).
\]
By \cite[(29)]{ACL}, two such characters agree if and only if their
shifted parameters belong to the same Weyl orbit.
The lowest conformal-weight space of $S_\lambda$ is $\CC u_\lambda$.
Zhu's correspondence therefore gives
$S_\lambda\cong S_\mu$ if and only if
$\lambda+\rho\in W(\mu+\rho)$.
\end{proof}

We give the complete orbit calculation for the set $\Lambda$,
using the orthonormal coordinates of Section~\ref{sec:affine-AV}.
For $q=3(\lambda+\rho)$, let $A(q)$ be the nondecreasing
rearrangement of $(|q_1|,|q_2|,|q_3|,|q_4|)$ and set
$P(q)=q_1q_2q_3q_4$.
Since $W(D_4)$ consists of signed permutations with an even
number of sign changes, the pair $(A(q),P(q))$ determines the
orbit of $q$.

\begin{theorem}\label{thm:seventeen-simple-quotients}
As $\lambda$ varies over all $405$ elements of $\Lambda$, the
$\mathcal W$-modules $S_\lambda$ give precisely $17$ isomorphism classes.
Representatives and the numbers of highest weights in $\Lambda$ giving
each class are listed in Table~\ref{tab:seventeen-simple-modules}.
\end{theorem}

\begin{proof}
Let $\lambda_1,\ldots,\lambda_{17}$ be the weights in the first
column of Table~\ref{tab:seventeen-simple-modules}, and put
\[
 \Lambda_i=\Lambda\cap\bigl(W(\lambda_i+\rho)-\rho\bigr).
\]
The explicit classification in \cite[Proposition~3.4 and Theorem~3.6]{AV}
gives $|\Lambda_i|=N_i$, where $N_i$ is the entry in the last column.
The pairs $(A(3(\lambda_i+\rho)),P(3(\lambda_i+\rho)))$ are distinct.
Thus the sets $\Lambda_i$ are disjoint, and
\[
 \sum_{i=1}^{17}|\Lambda_i|
 =1+3\cdot8+4+3\cdot8+24+3\cdot32+3\cdot24+96+64
 =405=|\Lambda|.
\]
It follows that
\[
 \Lambda=\coprod_{i=1}^{17}\Lambda_i.
\]
Proposition~\ref{prop:cyclic-coset-modules} now proves the assertion.
\end{proof}

\begin{table}[htbp]
\centering
\small
\renewcommand{\arraystretch}{1.15}
\begin{tabular}{c|c|r|c|r}
$\lambda$ & $A(3(\lambda+\rho))$ & $P(3(\lambda+\rho))$
 & $h_\lambda$ & $N$\\ \hline
$0$ & $(0,3,6,9)$ & $0$ & $0$ & $1$\\
$-\frac43\omega_1$ & $(0,3,5,6)$ & $0$ & $-1$ & $8$\\
$-\frac43\omega_3$ & $(1,2,4,7)$ & $56$ & $-1$ & $8$\\
$-\frac43\omega_4$ & $(1,2,4,7)$ & $-56$ & $-1$ & $8$\\
$-\omega_2$ & $(0,3,3,6)$ & $0$ & $-\frac97$ & $4$\\
$-\frac83\omega_1$ & $(0,1,3,6)$ & $0$ & $-\frac{10}7$ & $8$\\
$-\frac83\omega_3$ & $(1,2,4,5)$ & $-40$ & $-\frac{10}7$ & $8$\\
$-\frac83\omega_4$ & $(1,2,4,5)$ & $40$ & $-\frac{10}7$ & $8$\\
$-\frac43\omega_2$ & $(0,2,3,5)$ & $0$ & $-\frac{11}7$ & $24$\\
$-\frac43(\omega_1+\omega_3)$ & $(1,2,3,4)$ & $24$ & $-\frac{12}7$ & $32$\\
$-\frac43(\omega_1+\omega_4)$ & $(1,2,3,4)$ & $-24$ & $-\frac{12}7$ & $32$\\
$-\frac43(\omega_3+\omega_4)$ & $(0,1,2,5)$ & $0$ & $-\frac{12}7$ & $32$\\
$-\frac23\omega_1-\frac43\omega_2$ & $(0,2,3,3)$ & $0$ & $-\frac{13}7$ & $24$\\
$-\frac23\omega_3-\frac43\omega_2$ & $(1,1,2,4)$ & $8$ & $-\frac{13}7$ & $24$\\
$-\frac23\omega_4-\frac43\omega_2$ & $(1,1,2,4)$ & $-8$ & $-\frac{13}7$ & $24$\\
$-\frac73\omega_2$ & $(0,1,2,3)$ & $0$ & $-2$ & $96$\\
$-\frac23\omega_1-\frac53\omega_2-\frac23\omega_3$
 & $(0,1,1,2)$ & $0$ & $-\frac{15}7$ & $64$
\end{tabular}
\caption{The $17$ simple $\mathcal W$-modules obtained from highest-weight
vectors of $L_{-14/3}(\g)$-modules. The number $N$ counts weights in $\Lambda$ in
the indicated shifted Weyl orbit.}
\label{tab:seventeen-simple-modules}
\end{table}

The logarithmic $\mathcal W$-module and the seventeen simple
$\mathcal W$-modules give the following quotient of $A(\mathcal W)$.

\begin{theorem}\label{thm:A-W-quotient}
	There is a surjective algebra homomorphism
	\begin{equation}\label{eq:A-W-quotient}
		A(\mathcal W)\longrightarrow
		\mathbb C^{16}\times\mathbb C[\varepsilon]/(\varepsilon^2).
	\end{equation}
	In particular, $A(\mathcal W)$ is not semisimple and
	$\dim A(\mathcal W)\geq18$.
\end{theorem}
\begin{proof}
	Write $D=\mathbb C[\varepsilon]/(\varepsilon^2)$.
	The lowest generalized $L(0)$-eigenspace of the logarithmic $\mathcal W$-module
	$\mathcal C$ is the free rank-one $D$-module $D$.
	Its $A(\mathcal W)$-action commutes with multiplication by
	$\varepsilon$, and hence gives a homomorphism
	\[
	\phi:A(\mathcal W)\longrightarrow
	\operatorname{End}_D(D)=D.
	\]
	By Theorem~\ref{thm:log}, $L(0)=-2+\frac37\varepsilon$ on this space, so
	\[
	\phi([\omega]+2[\one])=\frac37\varepsilon.
	\]
	Thus $\phi$ is surjective. Its reduction modulo $\varepsilon$
	is the character corresponding to
	$\lambda_*=-\frac73\omega_2$, by the homomorphism
	$A(\mathcal W^{-38/7}(\g,f_{\mathrm{reg}}))\longrightarrow Z(U(\mathfrak{g}))$ constructed in the proof of
	Proposition~\ref{prop:cyclic-coset-modules}.
	
	Let $\mathfrak a=\ker\phi$, and let
	$\mathfrak m_1,\ldots,\mathfrak m_{16}$ be the kernels of the
	other sixteen characters in
	Theorem~\ref{thm:seventeen-simple-quotients}.
	The element $([\omega]+2[\one])^2$ belongs to $\mathfrak a$.
	On each of these sixteen characters it takes the value $(h+2)^2$,
	where $h$ is the corresponding lowest conformal weight.
	Table~\ref{tab:seventeen-simple-modules} gives $h\ne-2$ in every case.
	Thus $\mathfrak a+\mathfrak m_i=A(\mathcal W)$ for each $i$.
	The $\mathfrak m_i$ are distinct maximal ideals, hence pairwise
	comaximal. The Chinese remainder theorem gives \eqref{eq:A-W-quotient}.
\end{proof}

The bound $\dim R_{\mathcal W}\leq18$ in Proposition~\ref{prop:dim-RW-kvocijent} will imply that
\eqref{eq:A-W-quotient} is an isomorphism.

\section{\texorpdfstring{The Zhu algebra and the $C_2$-algebra of
 $\mathcal W$}{The Zhu algebra and the C2-algebra of W}}
\label{sec:singular-relations}

In this section we complete the proof of Theorem~\ref{thm:main}.
We determine the Zhu algebra $A(\mathcal W)$ and the Poisson algebra
$R_{\mathcal W}$, proving that
\[
 A(\mathcal W)\cong
 \mathbb C^{16}\times\mathbb C[\varepsilon]/(\varepsilon^2),
 \qquad
 \dim A(\mathcal W)=\dim R_{\mathcal W}=18.
\]
In particular, $\mathcal W$ is $C_2$-cofinite.
The description of $A(\mathcal W)$ also shows that the seventeen
simple $\mathcal W$-modules constructed in Section~\ref{sec:simple-modules}
exhaust the irreducible $\mathcal W$-modules.
Non-rationality follows from the existence of the  logarithmic $\mathcal W$-module
constructed in Section~\ref{sec:logarithmic}.

\subsection{\texorpdfstring{Singular vectors in $\mathcal W^{-38/7}(\g,f_{\mathrm{reg}})$ of conformal weight eight}{Singular vectors in W(-38/7)(D4, f reg) of conformal weight eight}}
Let $J$ be the maximal ideal of $\mathcal W^{-38/7}(\g,f_{\mathrm{reg}})$, so that
\[
 \mathcal W=\mathcal W^{-38/7}(\g,f_{\mathrm{reg}})/J.
\]
Write $J_n=J\cap \bigl(\mathcal W^{-38/7}(\g,f_{\mathrm{reg}})\bigr)_n$. We first locate three singular vectors of conformal weight eight
using principal Drinfeld--Sokolov reduction, and then give their explicit formula using the Miura homomorphism.

\begin{proposition}\label{prop:sing-v-existence}
	There exist linearly independent vectors $s_1,s_3,s_4\in J_8$
	such that, for each $i\in\{1,3,4\}$,  $\mathbb Cs_i$
	is the lowest conformal-weight space of a highest-weight
	$\mathcal W^{-38/7}(\g,f_{\mathrm{reg}})$-submodule of $\mathcal W^{-38/7}(\g,f_{\mathrm{reg}})$. The corresponding $A(\mathcal W^{-38/7}(\g,f_{\mathrm{reg}}))$-character is
	$\gamma_{2\omega_i-\frac47\rho}$, where $\gamma_\nu$ denotes
	the central character of the $\mathfrak g$-Verma module of highest weight $\nu$.
\end{proposition}
\begin{proof}
For the affine weights in this proof, adjoin the degree derivation $d$
with $[d,x(r)]=r x(r)$. Let $\Lambda_0,\delta$ vanish on $\h$ and
be dual to $K,d$, respectively. Extend finite weights by zero on $K,d$
and put $\widehat\rho=\rho+6\Lambda_0$.
By \cite[Theorem~1.2(5)]{JS}, the affine vertex algebra
	$V^{-38/7}(\g)$ contains singular vectors $v_i$, $i=1,3,4$,
	of weights
	\[
	\Lambda_i=-\frac{38}{7}\Lambda_0-14\delta+2\omega_i.
	\]
	Let $N_i$ be the $\widehat{\mathfrak g}$-submodule generated by $v_i$.
	It is a quotient of the Verma module $M(\Lambda_i)$ and has
	simple quotient $L(\Lambda_i)$.
	
	We verify the condition of \cite[(385)]{Ara07}. For every positive
	root $\beta$ of $\mathfrak g$ and $1\le n\le\operatorname{ht}\beta$,
	we have $n\le5$ and
	\begin{equation}
		\big\langle\Lambda_i+\widehat\rho,(-\beta+n\delta)^\vee\big\rangle
		=\frac{4n}{7}-\langle 2\omega_i+\rho,\beta^\vee\rangle
		\notin\mathbb Z.
	\end{equation}
	The same condition holds for $-\frac{38}{7}\Lambda_0$.
	By \cite[Theorem 9.1.3]{Ara07}, the functor $H^0_+$
	is exact on the corresponding blocks of category $\mathcal O$.
	We therefore have surjections
	\[
	H^0_+(M(\Lambda_i))\twoheadrightarrow H^0_+(N_i)
	\twoheadrightarrow H^0_+(L(\Lambda_i))
	\]
	and an embedding $H^0_+(N_i)\hookrightarrow \mathcal W^{-38/7}(\g,f_{\mathrm{reg}})$.
	By \cite[Propositions 9.2.1 and 9.2.3, Lemma 9.2.4]{Ara07},
	the first and last terms have one-dimensional lowest
	conformal-weight spaces, both with character
	$\gamma_{2\omega_i-\frac47\rho}$. Their lowest conformal weight is
	\[
	\frac78(2\omega_i\mid 2\omega_i+2\rho)
	-\langle 2\omega_i,\rho^\vee\rangle=14-6=8.
	\]
	It follows that $H^0_+(N_i)$ has no conformal weights below
	eight and has a one-dimensional, nonzero weight-eight space.
Choose $0\ne s_i\in H^0_+(N_i)_8\subset\W^{-38/7}(\g,f_{\reg})$.
	
		The weights $2\omega_i+\frac37\rho$, $i=1,3,4$, are distinct
and strictly dominant. Hence they belong to distinct Weyl orbits.
By the Harish--Chandra isomorphism, the central characters
$\gamma_{2\omega_i-\frac47\rho}$ are therefore pairwise distinct.
It follows that $s_1,s_3,s_4$ are linearly independent.
	
	Finally, let $\mathcal I$ be the maximal ideal of $V^{-38/7}(\g)$.
	Each $N_i$ is contained in $\mathcal I$.
	Since $-\frac{38j}{7}\notin\mathbb Z$ for $1\le j\le6$, \cite[Corollary 9.1.6(ii)]{Ara07} gives
	$H^0_+(L_{-38/7}(\g))\cong\mathcal W$.
	By exactness, the image of $H^0_+(\mathcal I)$ in $\mathcal W^{-38/7}(\g,f_{\mathrm{reg}})$
	is $J$. Hence $s_i\in J$.
\end{proof}

We express the singular vectors of Proposition~\ref{prop:sing-v-existence}
in terms of the primary generators $Q,P,R$ in \eqref{eq:primary-generators}.
Define the operator
\[
\begin{aligned}
 \mathcal D={}&88L(-4)+64L(-3)L(-1)-56L(-2)^2\\
 &+48L(-2)L(-1)^2-9L(-1)^4
\end{aligned}
\]
and the following vector of conformal weight eight:
\[
\begin{aligned}
 w={}&L(-2)R-\frac14L(-1)^2R\\
 &+\frac{16}{4747}\bigl(
 -12543054L(-8)+3623856L(-6)L(-2)\\
 &\qquad+2549580L(-5)L(-3)+1408005L(-4)^2\\
 &\qquad-825006L(-4)L(-2)^2-572268L(-3)^2L(-2)\\
 &\qquad+48287L(-2)^4\bigr)\mathbf1.
\end{aligned}
\]

\begin{proposition}\label{prop:explicit-singular-vectors}
The three singular vectors $s_1,s_3,s_4$ in
Proposition~\ref{prop:sing-v-existence} may be chosen as
\begin{equation}
\begin{aligned}
 s_1={}&P_{(-1)}P-2Q_{(-1)}Q
       +\frac53\mathcal D P+\frac13w,\\
 s_3={}&P_{(-1)}P+10P_{(-1)}Q-7Q_{(-1)}Q
       +\frac{10}{3}\mathcal D P-10\mathcal D Q-\frac43w,\\
 s_4={}&P_{(-1)}P-10P_{(-1)}Q-7Q_{(-1)}Q
       +\frac{10}{3}\mathcal D P+10\mathcal D Q-\frac43w.
\end{aligned}
\end{equation}
For $X\in\{P,Q,R\}$ and $n\in\ZZ$, write $X(n)=X_{(n+\wt(X)-1)}$.
The eigenvalues of $P(0)$, $Q(0)$ and $R(0)$ on the lines $\CC s_i$ are
	\begin{equation}\label{eq:eigenvalues}
		\begin{array}{c|ccc}
			&s_1&s_3&s_4\\ \hline
			P(0)&-4160&2080&2080\\
			Q(0)&0&-2080&2080\\
			R(0)&\dfrac{32065200128}{4747}&
			\dfrac{32065200128}{4747}&\dfrac{32065200128}{4747}
		\end{array}
	\end{equation}
With this notation the fields have expansions
	\[
 Y(P,z)=\sum_{n\in\ZZ}P(n)z^{-n-4},\ 
 Y(Q,z)=\sum_{n\in\ZZ}Q(n)z^{-n-4}, \ 
	Y(R,z)=\sum_{n\in\ZZ}R(n)z^{-n-6}.\]
\end{proposition}

\begin{proof}
For $a,b\in\{P,Q\}$ and $m=1,2$, the Virasoro commutator gives
\[
 L(m)(a_{(-1)}b)=(3m+4)a_{(m-1)}b.
\]
For instance, the expressions in \eqref{eq:primary-generators} and
the Virasoro relations give
\[
 P_{(1)}Q=224L(-2)Q-124L(-1)^2Q=\tfrac1{10}L(2)\mathcal DQ.
\]
Since $s_3-s_4=20(P_{(-1)}Q-\mathcal DQ)$, this gives
$L(2)(s_3-s_4)=0$. Applying the same commutator to the quadratic
terms and the Virasoro relations to $\mathcal DP,\mathcal DQ,w$ gives
\[
 L(1)s_i=L(2)s_i=0,\qquad i\in\{1,3,4\}.
\]
Thus $s_1,s_3,s_4$ are Virasoro-primary vectors of conformal
weight eight.
	
	Let
\[
 U_8=\operatorname{span}_{\mathbb C}
 \left\{L(-n_1)\cdots L(-n_r)X\ \middle|\
 \begin{gathered}
 X\in\{\mathbf1,P,Q,R\},\quad r\geq0,\\
 n_i\geq1\ (1\leq i\leq r),\quad
 \sum_{i=1}^r n_i=8-\mathrm{wt}(X)
 \end{gathered}\right\}.
\]
	The basis theorem for universal principal $W$-algebras
\cite[Proposition~4.12.1(iii)]{Ara07}, together with the changes of generators
established above, gives a basis of $\bigl(\mathcal W^{-38/7}(\g,f_{\mathrm{reg}})\bigr)_8$
consisting of the $22$ vectors in the first column of
Table~\ref{tab:weight-eight-pbw}.
Nineteen of these vectors belong to $U_8$, and the remaining three are
$P_{(-1)}P$, $P_{(-1)}Q$, and $Q_{(-1)}Q$.
Hence $\dim U_8=19$ and $\dim \bigl(\mathcal W^{-38/7}(\g,f_{\mathrm{reg}})\bigr)_8/U_8=3$.
The coefficients of $s_1,s_3,s_4$ relative to the classes of
$P_{(-1)}P$, $P_{(-1)}Q$, $Q_{(-1)}Q$ in this quotient form the matrix
\[
 \begin{pmatrix}
  1&0&-2\\
  1&10&-7\\
  1&-10&-7
 \end{pmatrix},
\]
whose determinant is $-100$. Therefore $s_1,s_3,s_4$ are linearly independent.
To apply the Kac determinant formula \cite{KacDet,FeFu}, write
$c=13-6(\tau+\tau^{-1})$. At $c=-50$, we have
$\tau+\tau^{-1}=21/2$ and $\tau$ is irrational. For positive integers
$r,s$, the weight
\[
 h_{r,s}=\frac{(r\tau-s)^2-(\tau-1)^2}{4\tau}
\]
is rational only when $r=s$, in which case $h_{r,r}=17(r^2-1)/8$.
Neither $4$ nor $6$ occurs. For highest weight zero only $(r,s)=(1,1)$
occurs, and the corresponding singular vector vanishes in the vacuum
module by $L(-1)\one=0$. Thus none of the Virasoro modules
generated by $\one,P,Q,R$ contains a nonzero primary vector of weight
eight. These modules form a direct sum by the basis theorem above.
It follows that the natural projection
\[
 \bigl(\mathcal W^{-38/7}(\g,f_{\mathrm{reg}})\bigr)_8
 \longrightarrow
 \bigl(\mathcal W^{-38/7}(\g,f_{\mathrm{reg}})\bigr)_8/U_8
\]
is injective on the space
	of Virasoro-primary vectors of weight eight,
	so this space has dimension at most three.
	Consequently, $s_1,s_3,s_4$ form a basis of that space.
	
	By Proposition \ref{prop:sing-v-existence}, there are three linearly independent singular vectors of conformal weight 8 in the maximal ideal $J$. Since every singular vector is
	Virasoro-primary, we conclude that $s_1,s_3,s_4$ are in $J$ and are annihilated by $P(n)$, $Q(n)$, $R(n)$, for all $n \geq 1$. Their eigenvalues in \eqref{eq:eigenvalues} are obtained by applying
	$P_{(3)}$, $Q_{(3)}$ and $R_{(5)}$ to these expressions in $\mathcal H$,
	using $[b_i(m),b_j(n)]=28m\delta_{ij}\delta_{m+n,0}\id$
	for $1\leq i,j\leq4$ and $m,n\in\ZZ$.
For $\eta_i=7(2\omega_i+\tfrac37\rho)$, the orthonormal coordinates of
$\eta_1,\eta_3,\eta_4$ are $(23,6,3,0)$, $(16,13,10,-7)$ and
$(16,13,10,7)$, respectively. Evaluation of $[P]$ and $[Q]$ at
these points in the shifted coordinates of Section~\ref{sec:Zhu-relations}
gives $(-4160,0)$, $(2080,-2080)$ and $(2080,2080)$.
Comparison with \eqref{eq:eigenvalues} identifies
$\mathbb C s_i$ with the highest-weight line of character
$\gamma_{2\omega_i-\frac47\rho}$ in Proposition~\ref{prop:sing-v-existence}.
\end{proof}

\begin{proposition}\label{prop:pbw-singular-vectors}
Table~\ref{tab:weight-eight-pbw} gives the coefficients of the singular
vectors $s_1,s_3,s_4$ of Proposition~\ref{prop:explicit-singular-vectors}
in the basis of conformal weight eight listed in its first column.
\end{proposition}
\begingroup
\small
\renewcommand{\arraystretch}{1.15}
\setlength{\LTleft}{\fill}
\setlength{\LTright}{\fill}

\begin{longtable}{@{}lrrr@{}}
	\caption{Expansions of the singular vectors
		$s_1,s_3,s_4$.}
\label{tab:weight-eight-pbw}\\
\toprule
Basis vector & $s_1$ & $s_3$ & $s_4$\\
\midrule
\endfirsthead

\toprule
Basis vector & $s_1$ & $s_3$ & $s_4$\\
\midrule
\endhead

\endfoot
\bottomrule
\endlastfoot

$P(-4)^2\mathbf1$
& $1$ & $1$ & $1$\\
$P(-4)Q(-4)\mathbf1$
& $0$ & $10$ & $-10$\\
$Q(-4)^2\mathbf1$
& $-2$ & $-7$ & $-7$\\
$L(-4)P(-4)\mathbf1$
& $\frac{440}{3}$ & $\frac{880}{3}$ & $\frac{880}{3}$\\
$L(-3)P(-5)\mathbf1$
& $\frac{320}{3}$ & $\frac{640}{3}$ & $\frac{640}{3}$\\
$L(-2)^2P(-4)\mathbf1$
& $-\frac{280}{3}$ & $-\frac{560}{3}$ & $-\frac{560}{3}$\\
$L(-2)P(-6)\mathbf1$
& $160$ & $320$ & $320$\\
$P(-8)\mathbf1$
& $-360$ & $-720$ & $-720$\\
$L(-4)Q(-4)\mathbf1$
& $0$ & $-880$ & $880$\\
$L(-3)Q(-5)\mathbf1$
& $0$ & $-640$ & $640$\\
$L(-2)^2Q(-4)\mathbf1$
& $0$ & $560$ & $-560$\\
$L(-2)Q(-6)\mathbf1$
& $0$ & $-960$ & $960$\\
$Q(-8)\mathbf1$
& $0$ & $2160$ & $-2160$\\
$L(-2)R(-6)\mathbf1$
& $\frac13$ & $-\frac43$ & $-\frac43$\\
$R(-8)\mathbf1$
& $-\frac16$ & $\frac23$ & $\frac23$\\
$L(-8)\mathbf1$
& $-\frac{66896288}{4747}$
& $\frac{267585152}{4747}$
& $\frac{267585152}{4747}$\\
$L(-6)L(-2)\mathbf1$
& $\frac{19327232}{4747}$
& $-\frac{77308928}{4747}$
& $-\frac{77308928}{4747}$\\
$L(-5)L(-3)\mathbf1$
& $\frac{13597760}{4747}$
& $-\frac{54391040}{4747}$
& $-\frac{54391040}{4747}$\\
$L(-4)^2\mathbf1$
& $\frac{7509360}{4747}$
& $-\frac{30037440}{4747}$
& $-\frac{30037440}{4747}$\\
$L(-4)L(-2)^2\mathbf1$
& $-\frac{4400032}{4747}$
& $\frac{17600128}{4747}$
& $\frac{17600128}{4747}$\\
$L(-3)^2L(-2)\mathbf1$
& $-\frac{3052096}{4747}$
& $\frac{12208384}{4747}$
& $\frac{12208384}{4747}$\\
$L(-2)^4\mathbf1$
& $\frac{772592}{14241}$
& $-\frac{3090368}{14241}$
& $-\frac{3090368}{14241}$\\

\end{longtable}
\endgroup

\begin{proof}
For $X\in\{P,Q,R\}$ and $n\geq0$, translation covariance gives
\[
 L(-1)^nX=n!\,X(-\wt(X)-n)\mathbf1.
\]
Substitution of the definitions of $\mathcal D$ and $w$
into $s_1,s_3,s_4$, followed by this identity, gives the
coefficients in Table~\ref{tab:weight-eight-pbw}.
\end{proof}

\subsection{\texorpdfstring{Structure of $A(\mathcal W)$ and $R_{\mathcal W}$}{Structure of A(W) and R(W)}}
\label{sec:Zhu-relations}

As in Section~\ref{subsec:principal-D4}, let $U=\W^{-38/7}(\g,f_{\reg})$.
In $A(\mathcal H)$ set
\[
 \xi_i=[b_i(-1)\one]-3\rho_i[\one]\quad(1\leq i\leq4),
 \qquad (\rho_1,\rho_2,\rho_3,\rho_4)=(3,2,1,0).
\]
Under the embedding $A(\mu):A(U)\hookrightarrow A(\mathcal H)$,
$[a]\mapsto[\mu(a)]$ \cite[Lemma~6.1]{ACL}, write $x,y,z,t$ for
the four polynomials in \eqref{eq:x-y-z-t} evaluated at
$(\xi_1,\xi_2,\xi_3,\xi_4)$. Their leading symbols are the generators
$x,y,z,t$ of $R_U$. The formulas for $\omega,P,Q$ give
\[
 [\omega]=\frac{x-126}{56},\qquad [P]=\frac{4y-x^2}{56},\qquad [Q]=\frac z7.
\]

\begin{proposition}\label{prop:classes-singular-vectors}
For the singular vectors $s_1,s_3,s_4\in J$ given in
Table~\ref{tab:weight-eight-pbw}, the following identities hold in $A(U)$:
	\begin{align*}
		12544[s_1]={}&128xt-5376t+15x^4-840x^3
		+15288x^2-79200x+299376\\
		&+64y^2-80x^2y+3136xy-38080y-512z^2,\\[2mm]
		-196[s_3]={}&8xt-336t-y^2+28z^2-14xy
		+1190y+63x^2-4950x+18711\\
		&-20yz+420xz-7140z,\\[2mm]
		-196[s_4]={}&8xt-336t-y^2+28z^2-14xy
		+1190y+63x^2-4950x+18711\\
		&+20yz-420xz+7140z.
	\end{align*}
	In particular, the three polynomials on the right-hand sides vanish
	in $A(\mathcal W)$.
\end{proposition}	
	
\begin{proof}
Substitution of \eqref{eq:primary-generators} in Table~\ref{tab:weight-eight-pbw}
gives $\mu(s_i)\in\mathcal H$ for $i=1,3,4$. In
$A(\mathcal H)=\CC[\xi_1,\xi_2,\xi_3,\xi_4]$ we have
	\[
	\bigl[
	b_{i_1}(-n_1)\cdots b_{i_r}(-n_r)\mathbf1
	\bigr]
	=
	(-1)^{\sum_{j=1}^r(n_j-1)}
    \prod_{j=1}^r(\xi_{i_j}+3\rho_{i_j}),
	\]
for $r\geq0$, $n_j\geq1$ and $1\leq i_j\leq4$ ($1\leq j\leq r$).
Applying this identity to $\mu(s_i)$ gives $A(\mu)([s_i])$. Expressing
these three invariant polynomials in $x,y,z,t$ gives the stated
formulas. Since $s_i\in J$, the quotient map $A(U)\to A(\mathcal W)$
sends each $[s_i]$ to zero.
\end{proof}

In $R_U=\CC[x,y,z,t]$, set
\begin{align}
 u_1&\coloneqq 12544\,\bar s_1
     =128xt+15x^4-80x^2y+64y^2-512z^2,\label{eq:u-one}\\[2pt]
 u_3&\coloneqq -196\,\bar s_3=8xt-y^2+28z^2-20yz,\label{eq:u-three}\\[2pt]
 u_4&\coloneqq -196\,\bar s_4=8xt-y^2+28z^2+20yz.\label{eq:u-four}
\end{align}
The polynomial expressions are the weight-eight components of the
corresponding polynomials in Proposition~\ref{prop:classes-singular-vectors}.
Let $I=\ker(R_U\twoheadrightarrow R_{\mathcal W})=(J+C_2(U))/C_2(U)$.
Since $s_1,s_3,s_4\in J$, we have $u_1,u_3,u_4\in I$.
	
All three polynomials vanish on the $t$-axis. Hence $(u_1,u_3,u_4)$
has infinite codimension in $R_U$.
Arakawa, Lam and Yamada used the mode $(W^3)_{(1)}$ to obtain
relations in the $C_2$-algebra of the $\mathfrak{sl}_2$ parafermion
algebra \cite[Proposition~5.4 and Section~6]{ALY}.

For $a\in U$ and $b\in J$, \eqref{eq:first-product-miura} gives
$28B_0(\bar a,\bar b)=\overline{a_{(1)}b}\in I$ since $a_{(1)}b\in J$.
Thus $B_0(R_U,I)\subseteq I$.

\begin{proposition}\label{prop:further-relations}
The ideal $I$ contains the polynomials
\[
\begin{aligned}
 v_1&\coloneqq 10xyz+12zt,\\
 v_2&\coloneqq yt+6xz^2,\\
 v_3&\coloneqq 96t^2+16xyt+96x^2z^2+192yz^2,\\
 v_4&\coloneqq 32x^3y-128xy^2-192x^2t+384yt-2304xz^2.
\end{aligned}
\]
\end{proposition}

\begin{proof}
By \eqref{eq:u-one}--\eqref{eq:u-four},
\[
\begin{aligned}
 \tfrac1{40}(u_4-u_3)
  &=\tfrac{49}{10}\overline{s_3-s_4}=yz,\\
 \tfrac12(u_3+u_4)
  &=-98\overline{s_3+s_4}=8xt-y^2+28z^2,\\
 \tfrac15(u_1-8u_3-8u_4)
  &=\tfrac{1568}{5}\overline{8s_1+s_3+s_4}
    =3x^4-16x^2y+16y^2-192z^2.
\end{aligned}
\]
Equation~\eqref{eq:x-y-z-t} gives
\[
 B_0(z,y)=2z\sum_{i=1}^4(x-q_i^2)=6xz,\qquad
 B_0(z,z)=\sum_{i=1}^4\prod_{j\ne i}q_j^2=t.
\]
Since $yz \in I$ and $B_0(R_U,I)\subseteq I$, we have
$B_0(z,yz)\in I$. The Leibniz rule gives
\[
 B_0(z,yz)=B_0(z,y)z+yB_0(z,z)=6xz^2+yt=v_2.
\]
For the other three polynomials, \eqref{eq:x-y-z-t} gives
\[
\begin{array}{c|cccc}
 B_0 & x & y & z & t\\ \hline
 y & 8y & 4xy+12t & 6xz & 8xt+16z^2\\
 t & 12t & 8xt+16z^2 & 4yz & 4yt+12xz^2
\end{array}
\]
and hence
\[
 v_1=B_0(y,yz),\qquad
 v_3=\tfrac12 B_0(t,u_3+u_4),\qquad
 v_4=\tfrac15 B_0(y,u_1-8u_3-8u_4).
\]
Since $yz,u_3+u_4,u_1-8u_3-8u_4\in I$ and $B_0(R_U,I)\subseteq I$,
all four polynomials belong to $I$.
\end{proof}

\begin{proposition}\label{prop:dim-RW-kvocijent}
Let $I_0=(u_1,u_3,u_4,v_1,v_2,v_3,v_4)\subset R_U$.
The classes of the eighteen monomials
\[
 1,x,\ldots,x^7,\qquad
 y,xy,x^2y,x^3y,y^2,\qquad
 z,xz,x^2z,x^3z,\qquad t
\]
span $R_U/I_0$. In particular, $\dim R_{\mathcal W}\leq18$.
\end{proposition}

\begin{proof}
In $R_U/I_0$, the relations $u_1=u_3=u_4=0$ give
\[
 yz=0,\qquad
 xt=-\frac16y^2+\frac7{24}x^2y-\frac7{128}x^4,\qquad
 z^2=\frac1{12}y^2-\frac1{12}x^2y+\frac1{64}x^4.
\]
Substitution in $v_1=v_2=v_4=0$ gives
\[
 zt=0,\qquad
 yt=\frac18x^3y-\frac{19}{640}x^5,\qquad
 xy^2=\frac34x^3y-\frac{41}{320}x^5.
\]
Multiplication of the formula for $xt$ by $z$ gives $x^4z=0$.
Multiplication of the formula for $z^2$ by $y$ gives
\[
 y^3=x^2y^2-\frac3{16}x^4y.
\]
Using these identities in $y(xt)-x(yt)$ gives
\[
 0=\frac7{512}(x^6-4x^4y).
\]
Thus $x^4y=x^6/4$, and substitution in the preceding formulas yields
\[
 y^3=\frac1{80}x^6,\qquad
 xyt=\frac1{640}x^6,\qquad
 x^2z^2=-\frac1{3840}x^6.
\]
It follows that $v_3=96t^2$ in $R_U/I_0$, so $t^2=0$.
Also $x^8=-3840x^4z^2=0$ since $x^4z=0$.

The relations for $t^2,yt,zt,yz,z^2,y^3$ show that $R_U/I_0$ is
spanned over $\CC[x]$ by $1,y,y^2,z,t$. The formulas for $xt$ and
$xy^2$ reduce the multiples of $t$ and $y^2$ by positive powers
of $x$. Finally, $x^8=0$, $x^4y=x^6/4$ and $x^4z=0$ give the
eighteen monomials in the statement. Since $I_0\subseteq I$, the
surjection $R_U/I_0\twoheadrightarrow R_{\mathcal W}$ gives
$\dim R_{\mathcal W}\leq18$.
\end{proof}

Proposition~\ref{prop:dim-RW-kvocijent} gives the following consequence.

\begin{corollary}
	The simple principal $W$-algebra
	$\mathcal W$ is $C_2$-cofinite.
\end{corollary}

\begin{theorem}\label{thm:Zhu-C2-iso}
	We have
	\[
	A(\mathcal W)\cong
	\mathbb C^{16}\times
	\mathbb C[\varepsilon]/(\varepsilon^2).
	\]
	In particular, $\dim A(\mathcal W)=18$ and the seventeen $\mathcal W$-modules
	listed in Table~\ref{tab:seventeen-simple-modules} exhaust the irreducible
$\mathcal W$-modules. The $C_2$-algebra of $\mathcal W$ has the presentation
	\[
	R_{\mathcal W}\cong
	\frac{\mathbb C[x,y,z,t]}
	{(u_1,u_3,u_4,v_1,v_2,v_3,v_4)}.
	\]
Here $u_1,u_3,u_4$ are defined in \eqref{eq:u-one}--\eqref{eq:u-four}
and $v_1,v_2,v_3,v_4$ in Proposition~\ref{prop:further-relations}.
In particular, $\dim R_{\mathcal W}=18$, and its
		Hilbert series is
		\[
		\operatorname{Hilb}_{R_{\mathcal W}}(q)
		=
		1+q^2+3q^4+4q^6+4q^8
		+3q^{10}+q^{12}+q^{14}.
		\]
	\end{theorem}
	
\begin{proof}
	By Theorem \ref{thm:A-W-quotient} and Proposition \ref{prop:dim-RW-kvocijent}, we have
	$	18\leq \dim A(\mathcal W) \leq\dim R_{\mathcal W}\leq18$.
	Thus $\dim A(\mathcal W)=\dim R_{\mathcal W}=18$, and the surjection
	in Theorem \ref{thm:A-W-quotient} is an isomorphism, giving
	$	A(\mathcal W)\cong\mathbb C^{16}\times
	\mathbb C[\varepsilon]/(\varepsilon^2).	$
	Moreover, since $I_0 \subseteq I$,
	the surjection $\mathbb C[x,y,z,t]/I_0\to R_{\mathcal W}$ is also an
	isomorphism by Proposition \ref{prop:dim-RW-kvocijent}.
	The monomials listed there therefore form a basis
	of $R_{\mathcal W}$ and give the asserted Hilbert series.
	
	Finally, $A(\mathcal W)$ has exactly $17$ isomorphism classes of irreducible $A(\mathcal W)$-modules.
	Since $\mathcal W$ is $C_2$-cofinite, every irreducible weak
	$\mathcal W$-module is ordinary by \cite[Corollary~5.8]{ABD}.
	By Zhu's correspondence \cite[Theorem~2.2.2]{Zhu}, the $\mathcal W$-modules constructed
	in Theorem \ref{thm:seventeen-simple-quotients} exhaust the irreducible
	$\mathcal W$-modules.
\end{proof}

\subsection*{Acknowledgements.}
The authors are partially supported by the Croatian Science Foundation
under the project IP-2022-10-9006 and by the project
``Implementation of cutting-edge research and its application as part
of the Scientific Center of Excellence for Quantum and Complex Systems,
and Representations of Lie Algebras'', Grant No.~\mbox{PK.1.1.10.0004},
co-financed by the European Union through the European Regional
Development Fund -- Competitiveness and Cohesion Programme 2021--2027.

D.A. was also partially  supported in part
by the European Union - NextGenerationEU through the National Recovery and Resilience
Plan 2021–2026, via an institutional grant of the University of Zagreb Faculty of Science, IK IA 1.1.3, Impact4Math.

I.V. was also partially supported by the European Union (NextGenerationEU) under the Croatian Recovery and Resilience Plan 2021–2026 (NRRP), through the University of Split institutional project “Mathematical modelling and simulations of physical systems (MaMoS) IP-UNIST-46”,
approved by the Ministry of Science, Education and Youth of the Republic
of Croatia.

 \subsection*{Disclosure of AI use.}
ChatGPT was used for mathematical discussions and assistance
with mathematical arguments, computations, and the preparation
of this manuscript.
The authors take full responsibility for the mathematical content.

\end{document}